\documentclass{elsarticle}
\usepackage{amsmath} %align
\usepackage{hhline}
\usepackage{natbib}
\usepackage{lineno, hyperref}
\usepackage{wrapfig}
\usepackage{bm,upgreek}
\usepackage{mathrsfs}
\usepackage{tablefootnote}
\usepackage{amsfonts}
\usepackage{color}
\usepackage{ulem}
\usepackage{accents}
\usepackage[utf8]{inputenc}
\usepackage{graphicx} %figure
\usepackage{float} %[H]
\usepackage{amssymb} %mathbb
\usepackage{subcaption}
\usepackage{algorithm2e} %algorithm
\usepackage{cleveref} %cref

\usepackage[export]{adjustbox}
\graphicspath{ {./images/} }
\RestyleAlgo{ruled} %algorithm style
\SetKwComment{Comment}{/* }{ */} %comment style in algorithm
\newcommand{\dx}{\mathrm{d}x}
\newcommand{\rx}{\mathrm{x}}
\newcommand{\rs}{\mathrm{s}}
\newcommand{\V}{\textbf{v}}
\newcommand{\vecko}{\textbf{v}}
\newcommand{\ucko}{\textbf{u}}
\newcommand{\pder}[2]{\frac{\partial #1}{\partial #2}}
\newcommand{\PHI}{\bm{\varphi}}
\newcommand{\U}{\textbf{u}}
\newcommand{\F}{\textbf{f}}
\newcommand{\N}{\textbf{n}}

\newcommand{\Div}[1]{\text{div}\left(#1\right)}
\renewcommand{\div}{\text{div}}
\newcommand{\tr}{\text{tr}}
\newcommand{\FF}{\mathbb{F}}
\newcommand{\DD}{\mathbb{D}}
\newcommand{\BB}{\mathbb{B}}
\newcommand{\fid}[1]{\ensuremath{\accentset{\triangledown}{#1}}}
\newcommand{\oldB}{\fid{\BB}}
\newcommand{\LL}{\mathbb{L}}
\newcommand{\II}{\mathbb{I}}
\newcommand{\TT}{\mathbb{T}}
\newcommand{\PP}{\mathbb{P}}
\newcommand{\mesh}{\text{mesh}}
\newcommand{\norm}[1]{\lVert #1 \rVert}
\newcommand{\mycomment}[1]{}
\newcommand{\meshlev}[2]{\text{mesh}_{#1}^{#2}}
\newcommand{\noref}{\text{None}}
\newcommand{\qualityref}{\text{Quality}}
\newcommand{\eicref}{\text{Eikonal}}

\makeatletter
\renewenvironment{cases}[1][l]{\matrix@check\cases\env@cases{#1}}{\endarray\right.}
\def\env@cases#1{%
  \let\@ifnextchar\new@ifnextchar
  \left\lbrace\def\arraystretch{1.2}%
  \array{@{}#1@{\quad}l@{}}}
\makeatother

\begin{document}
\begin{frontmatter}
\title{Geometric re-meshing strategies to simulate contactless rebounds of elastic solids in fluids}
\author[FMP]{J. Fara}
\ead{fara@karlin.mff.cuni.cz}

\author[FMP,UU]{S. Schwarzacher}
\ead{schwarz@karlin.mff.cuni.cz}

\author[FMP]{K. T\r{u}ma\corref{cor1}}
\ead{ktuma@karlin.mff.cuni.cz}

\cortext[cor1]{Corresponding author. Tel.: (+420) 221 913 248.}

\address[FMP]{Faculty of Mathematics and Physics, Charles University, Sokolovsk\'{a} 83, 186 75 Prague, Czech Republic}
\address[UU]{Department of Mathematics, Box 480, 751 06 Uppsala, Sweden}

\begin{abstract}
  The paper deals with the rebound of an elastic solid off a rigid wall of a container filled with an incompressible Newtonian fluid. Our study focuses on a collision-free bounce, meaning a rebound without topological contact between the elastic solid and the wall. This has the advantage of omitting any artificial bouncing law. In order to capture the contact-free rebound for very small viscosities an adaptive numerical scheme is introduced. 

  The here-introduced scheme is based on a Glowinski time scheme and a localized arbitrary Lagrangian-Eulerian map on finite elements in space.   
  The absence of topological contact requires that very thin liquid channels are solved with sufficient accuracy. It is achieved via newly developed geometrically driven adaptive strategies. Using the numerical scheme, we present here a collection of numerical experiments. A rebound is simulated in the absence of topological contacts. Its physical relevance is demonstrated as, with decreasing viscosities, a free rebound in a vacuum is approached.
  Further, we compare the dynamics with a second numerical scheme; a here-introduced adaptive purely Eulerian level-set method. The scheme produced the same dynamics for large viscosities. However, as it requires a much higher computational cost, small viscosities can not be reached by this method.
  The experiments allow for a better understanding of the effect of fluids on the dynamics of elastic objects.
  Several observations are discussed, such as the amount of elastic and/or kinetic energy loss or the precise connection between the fluid pressure and the rebound of the solid.
\end{abstract}

\begin{keyword}
Contact \sep Numerical Simulation \sep Finite Element Method \sep Fluid-Structure Interaction \sep Rebound 
\end{keyword}

\end{frontmatter} 

\section{Introduction}

\noindent
This work is dedicated to the simulation of an elastic object's rebound off a container's rigid wall filled with an incompressible Newtonian fluid. The effect of the fluid on an elastic rebound is still not a fully described phenomenon despite its evident significance for so many important applications.
The most prominent strategy for such simulations is to
model the contact between solids inside fluids with additional heuristic constraints. This is a much-debated subject and not a straightforward task at all. Nevertheless, many recent numerical studies (\cite{zonca2021polygonal,Ager2019,Burman2022, frei2018FSIbook, richter}) could simulate the contact between two materials by artificially enforced contact laws. One example is the Nitsche method, which prevents materials from penetrating (\cite{richter,Ager2021,burman}).

{\em In this paper, we follow a different approach that avoids any additional constraints.}
We can do that as we base our method on analytical results \cite{MR2592281,gravina2022,GraHil16}, which show that no contact between a smooth solid and a rigid wall is possible in finite time.  

The contactless behavior of solids in fluids has been known for some time.
For the Stokes fluid, analytical estimates based on approximation techniques performed by Brenner
\cite{Brenner1961} and Cooley and O'Neill \cite{cooley_o-neill-1969} excluded the collision due to
the singular pressure appearing when the incompressible fluid is squeezed together.
Later, Hesla \cite{Hesla} and Hillairet \cite{MR2354496} showed that contact for smooth rigid objects is also impossible
for the incompressible Navier-Stokes fluid. See also \cite{gravina2022} for some analytical investigations on contactless rebounds. Nevertheless, mathematically, the question of how smooth elastic objects can collide (or bounce without contact) inside incompressible viscous fluids is still widely open. Hence numerical investigations are besides others of relevance for pure analytic questions.

Numerically, the contactless rebound was first observed already in 2016 by Frei \cite{frei2016eulerian} and later published in Chapter 3 of \cite{frei2018FSIbook}. After that, a first publicly available source code called LocModFE for interface problems dealing also with the contactless rebound appeared \cite{freirichterwick2021}.

{\em We aim to further confirm this result by numerical simulations and show that contactless rebounds are a valid numerically reasonable approximation strategy for fluid-structure interaction.} It comes, however, with the challenge of capturing the small fluid-channel between the solids, in particular, in the relevant case of small viscosities. In this paper, we introduce an adaptive and geometrical re-meshing strategy. This is implemented into an {\em Updated arbitrary Lagrangian-Eulerian (ALE) scheme.
The success of the strategy is demonstrated by numerical experiments that allow for viscosities as small as $10^{-3}$ Pa s. It is shown that with decreasing viscosities, the free bounce in a vacuum is indeed approached. The simulations further allow for a better understanding of the physical process.}

Various prominent FSI schemes have been tested on contact and bouncing problems. This, on the one hand, includes the fully Eulerian approach \cite{dunne2006eulerian,dunne2006adaptive,cottet2006eulerian,richter2010finite,richter2013fully,frei2016eulerian,frei2016long,gravina2022} and, on the other hand, the (Updated) ALE method \cite{hughes1981, Donea1982, stein2003,bungartz2006, HronTurek2006a, richter2010finite, HronTurekMadlikEtAl2010, turek2011numerical, frei2016long}.

In \cite{richter2013fully,richter} the authors do compare Eulerian schemes with ALE schemes and it seems that their conclusion is that the ALE scheme is not capable of capturing a large deformation close to the boundary and the simulations crash before contact of elastic solids in fluids is established. In this paper, we show that an Updated ALE method incorporated with an adaptive re-meshing strategy is indeed sufficiently accurate to obtain a contactless rebound. Actually, it allows for smaller viscosities than those computed in aformentioned papers.

Our leading experimental set up is introduced in Subsection~\ref{sec41}. It is a fluid-structure interaction (FSI) problem for a contactless rebound in a viscous incompressible fluid. The problem is designed to be as simple as possible to see the rebound undisturbed by other effects.

The no-contact pathology of the bounce in an incompressible fluid has important implications for the
numerical description of this problem. 
The fact that very narrow fluid channels appear
between the body and the fixed wall complicates the numerical solution. However, two important
advantages stem from this effect. First, technically, it allows for simulation without topological
changes in the fluid domain. Second, it allows for a simulation of an elastic rebound without
further ad-hoc contact laws or similar: The momentum equation, mass conservation, and the no-slip boundary conditions alone produce a proper rebound. One main motivation for the research performed here is {\em the great advantage of contactless strategies for computer simulations}. 
Indeed, excluding topological changes a priori excludes a common reason a computer code breaks down.

{\em The central technical achievement of the present work is a finite element ALE based adaptive program that allows for accurate simulation of this event for very small viscosities and respectively degenerate geometries and singular fluid pressures}.
For that, two strategies have been developed/implemented:

\begin{enumerate}
\item A macroscopic re-meshing strategy that adaptively modifies the reference configuration for the coupled system. This means that the ALE map is {\em localized in time}, with the advantage of remaining stable meshes even when the fluid domain changes drastically, as when a contact is approaching. 

\item A geometrically induced microscopic adaptive refinement that resolves the narrow liquid layer. In particular, a marker that is sensitive to the distance between solid objects.
\end{enumerate}

\noindent
Geometric re-meshing strategies seem to be new in this context. A respective library, ADmesh \cite{ADmesh}, compatible with FEniCS finite element code, has been created by one of the authors. The adaptive strategy is based on local mesh operations. See \cite{dolejvsi1998anisotropic, li20053d, de1999parallel} for further references on related geometric mesh adaptations.

The here-developed macroscopic/microscopic adaptive ALE scheme is very suitable for computing contactless rebounds. Its success is best demonstrated by its ability to allow for rather small viscosities compared to previous numerical experiments \cite{richter2013fully,frei2018FSIbook,richter}. Indeed, the viscosity can be as small as $10^{-3}$~Pa\,s with otherwise unchanged parameters. 
In this case, the distance between the elastic ball and the wall during the rebound becomes tens of micrometers, which highlights the necessity of the novel adaptive refinement method. 

Further, the macroscopic/microscopic adaptive ALE scheme's experiments reveal several observations with possible implications for physically understanding elastic objects surrounded by fluids. This includes the quantification of the loss of elastic energy over the rebound, the coincidence between pressure peeks and the elastic deformation, the dynamic of the geometry of the solid surface (its oscillation and/or concavities) when contact is approaching, and their relation to decreasing viscosity. For comparison, we also compute the rebound of the elastic ball from the wall in a vacuum. We observe that the results converge to the vacuum case as the viscosity decreases. 
Even so, this strongly emphasizes that the simulations are indeed reflecting reality. The precise regime of validity for contactless rebounds is still a widely open question, including how materials touch inside fluids. Accurate computer simulations seem to be crucial to further understanding its precise meaning and the resulting dynamical peculiarities. In particular, in the future, we aim for parameters where a simulation can withstand experimental comparison.

Finally, we follow the tradition of comparing Eulerian to ALE strategies. This we do on the one hand to confirm that the computed solutions are in accordance and correct. On the other hand to further promote the advantages of the here introduced macroscopic/mi\-cro\-sco\-pic adaptive ALE scheme. Our Eulerian scheme is based on the level-set method.

The here presented Eulerian scheme originates in \cite{gravina2022} but underwent a serious revision. Most importantly, we introduced an adaptive refinement around the smooth interface. For not too small viscosities it provides a very accurate simulation of the contactless rebound, which is in excellent accordance with the macroscopic/microscopic adaptive ALE scheme. However, due to its significantly higher computational effort, it cannot simulate the physically relevant case of small viscosities.

The paper is structured as follows. The next section deals with the modeling part, where both Lagrangian and Eulerian settings are introduced, as well as the constitutive relations for fluid and solid. Section~3 compares different ALE approaches, in this paper we employ the so-called {\em Updated ALE scheme} with {\em adaptive meshing strategy} based on on the meshing library ADmesh \cite{ADmesh}, which was developed by the authors for the scheme. This method is then applied to a FSI problem of contactless rebound of a ball in an incompressible Newtonian fluid. The experiment is described in Subsection 4.1. The rest of Section~4 is dedicated to the presentation of the numerical results of the method. Specifically, the convergence in space and time. Section~5 provides some physical aspects of the obtained results, specifically, the problem is computed with different values of material parameters. In particular, we investigate the behavior of decreasing viscosities (up to $10^{-3}$ Pa s) and show that it approaches the bounce obtained for a ball in a vacuum.
This follows Section~6, where we compare the ALE experiments with the experiments of the introduced level-set scheme. Concluding remarks can be found in Section~7.

\section{Modeling of rebound of the elastic ball in viscous fluid}
In this section, we formulate the model to describe the interaction between the elastic solid and the viscous fluid. Since the fluid is described naturally in the current (Eulerian) setting and the solid in the reference (Lagrangian) setting, we first recall the Eulerian-Lagrangian formulation. Next, we reformulate everything in the purely Eulerian setting. 

Both Eulerian and Lagrangian frameworks were introduced more than 200 years ago and once established, it was known how to transform between them. Theoretical works on large deformations of solids have been done in the Eulerian setting, see for example \cite{truesdell1965} where the solid is described both in the classical Lagrangian as well as Eulerian configurations. However, from the numerical point of view the solid has usually been studied in the reference configuration because it enables one to solve the problem on a fixed mesh. In the context of numerical analysis, the idea of fully Eulerian FSI first appeared in 2006 using the finite element method \cite{dunne2006eulerian,dunne2006adaptive} as well as the finite differences \cite{cottet2006eulerian}. Later, it was used in many numerical studies, see for instance \cite{richter2010finite,richter2013fully,frei2016eulerian,frei2016long}.

\subsection{Eulerian-Lagrangian fluid-structure interaction formulation}
We consider an incompressible Newtonian fluid occupying the domain $\Omega_{\rm f}(t)$ and the neo-Hookean solid occupying $\Omega_{\rm s}(t)$ such that both bodies fill the whole domain $\Omega=\Omega_{\rm f}(t)\cup\Omega_{\rm s}(t)$, where $\Omega\subset\mathbb{R}^d, d=2,3$. The homogeneous incompressible fluid is described by the balance of mass and linear momentum in the current configuration, i.e.,
\begin{align}
  \div_x\vecko&=0,\\
  \rho_{\rm f}\left(\pder{\vecko}{t}+\vecko\cdot\nabla_x\vecko\right)&=\div_x\TT_{\rm f},
\end{align}
where $\vecko$ is the fluid velocity, $\rho_{\rm f}$ a  constant fluid density, and $\TT_{\rm f}$ is the Cauchy stress tensor and $x\in\Omega_{\rm f}(t)$ is the position in the current configuration and all spatial derivatives are taken with respect to it. In the case of the incompressible Newtonian fluid, the Cauchy stress tensor reads
\begin{equation}\label{CauchyNewtonian}
  \TT_{\rm f}=-p_{\rm f}\II+2\mu\DD,
\end{equation}
where $p_{\rm f}$ is the unknown pressure, $\mu$ is a constant shear viscosity and $\DD=(\nabla_x\vecko+(\nabla_x\vecko)^{\rm T})$ is the symmetric part of the velocity gradient.

The solid is described by the balance of mass (that in the incompressible case reduces to a constant density) and linear momentum as well, but in its reference configuration, i.e.,
\begin{align}
  \rho_{\rm s}&=\frac{\rho_{\rm s, 0}}{J},\\
  \rho_{\rm s}\pder{^2\ucko}{t^2}&=\div_X\PP,
\end{align}
where $\rho_{\rm s}$ and $\rho_{\rm s,0}$ are the densities in the current and the reference configurations, respectively, $\ucko$ is the displacement of the solid and $J=\det(\FF)=\det(\II+\nabla_X \ucko)$ the determinant of the deformation gradient $\FF$, which equals one in the case of an incompressible solid. Finally, $\PP$ is the first Piola-Kirchhoff stress tensor, and $X\in\Omega_{\rm s}(t=0)$ is the position in the reference configuration.
In the case of compressible neo-Hookean solid, the first Piola-Kirchhoff stress is computed from the prescribed strain energy $W$ by
\begin{equation}
  \PP=\pder{W}{\FF},\quad W=\frac{G}{2}\left(\tr(\BB-\II)-\ln\det\BB\right)+\frac{\kappa}{2}(J-1)^2,
\end{equation}
where $\BB = \FF\FF^{\rm T}$ is the left Cauchy-Green tensor, $G$ is the elastic shear modulus and $\kappa$ is the elastic bulk modulus.
This results in the first Piola-Kirchhoff stress in the form
\begin{equation}
  \label{eq:comp}
  \PP = G (\FF - \FF^{-\rm T}) + \kappa (J - 1) J \FF^{-\rm T}.
\end{equation}

In the case of incompressible neo-Hookean solid, the first Piola-Kirchhoff stress is computed from the prescribed strain energy $W$ by
\begin{equation}
  \PP=\pder{L}{\FF},\quad L = W + P(J-1),\quad W=\frac{G}{2}\tr(\BB-\II),\label{firstPiola1}
\end{equation}
where $L$ is the Lagrange function that takes care of the incompressibility restriction $J=1$ using the Lagrange multiplier $P$.
This leads to
\begin{equation}
  \PP=G \FF + JP\FF^{-\rm T}.\label{firstPiola2}
\end{equation}
The compressible neo-Hookean solid \eqref{eq:comp} approximates the incompressible neo-Hookean solid \eqref{firstPiola2} for the elastic bulk modulus $\kappa$ large enough.

Finally, we have to prescribe the conditions on the interface $\Gamma_{\mathrm{int}}=\overline{\Omega_{\rm f}}\cap\overline{\Omega_{\rm s}}$, which consists of the kinematic interface condition (equality of velocities)
\begin{equation}
  \vecko = \pder{\U}{t}\quad{\textrm{on}}\ \Gamma_{\mathrm{int}}
\end{equation}
and the dynamic interface condition (equality of the tractions)
\begin{equation}
  \TT_{\rm f}{\bf n} = \PP  {\bf N}\quad{\text{on}}\ \Gamma_{\mathrm{int}},
\end{equation}
where ${\bf n}$ and ${\bf N}$ are the unit normals to $\Gamma$ in the current and reference configurations.

\subsection{Purely Eulerian FSI formulation}
Since the fluid is naturally given in a different (Eulerian) configuration than the solid (Lagrangian), it may be convenient to formulate both bodies in the same (Eulerian) framework. Thus, we reformulate the solid into the Eulerian setting where the balance of mass and linear momentum for an incompressible solid read
\begin{align}
  \div_x\vecko&=0,\label{balmasssolid}\\
  \rho_{\rm s}\left(\pder{\vecko}{t}+\vecko\cdot\nabla_x\vecko\right)&=\div_x\TT_{\rm s},\label{balmomsolid}
\end{align}
where $\TT_{\rm s}$ is the Cauchy stress tensor that is related to the first Piola-Kirchhoff stress tensor through
\begin{equation}
  \TT_{\rm s} = \frac{1}{J}\PP\FF^{\rm T} = P\II+G\BB = -p\II +G\BB^d\label{CauchySolid}
\end{equation}
where we substituted $\PP$ from \eqref{firstPiola2} and defined a new unknown pressure $p:=-P-\frac{G}{d}\tr\BB$. The symbol $\BB^d=\BB-\frac{1}{d}(\tr\BB)\II$ denotes the deviatoric part of the tensor. To close the system of equations \eqref{balmasssolid}--\eqref{CauchySolid}, we need to provide the equation for the left Cauchy-Green tensor $\BB$. We take the material time derivative (denoted by a dot) of the left Cauchy-Green tensor
\begin{equation*}
  \pder{\BB}{t}+\vecko\cdot\nabla\BB=\frac{{\rm d}\BB}{{\rm d}t}=:\dot{\BB}=\dot{\overline{\FF\FF^{\rm T}}}=\dot{\FF}\FF^{\rm T}+\FF\dot{\FF}^{\rm T}=
  \LL\FF\FF^{\rm T}+\FF\FF^{\rm T}\LL^{\rm T}=\LL\BB+\BB\LL^{\rm T},
\end{equation*}
where we used the relation between the rate of deformation gradient $\FF$ and velocity gradient $\LL:=\nabla\vecko$
\begin{equation}
  \dot{\FF}=\LL\FF.
\end{equation}
Thus, we obtained an evolution equation for the left Cauchy-Green tensor in the form
\begin{equation}
  \oldB:=\pder{\BB}{t}+\vecko\cdot\nabla\BB-\LL\BB-\BB\LL^{\rm T}=\mathbb{O}.
\end{equation}
Here, $\oldB$ is an objective upper convected Oldroyd derivative. For details, see \cite{Oldroyd1950}.
In summary, we obtained the governing equations that describe the evolution of incompressible neo-Hookean solid in the Eulerian setting
\begin{align}
  \div_x\vecko &= 0,\\
  \rho_{\rm s}\left(\pder{\vecko}{t}+\vecko\cdot\nabla_x\vecko\right)&=\div_x\TT_{\rm s},\quad \TT_{\rm s} = -p\II +G\BB^d,\\
  \pder{\BB}{t}+\vecko\cdot\nabla\BB-\LL\BB-\BB\LL^{\rm T}&=\mathbb{O}.
\end{align}
Since both fluid and solid are now formulated in the same Eulerian setting, the purely Eulerian FSI model can be written in the following way
\begin{align}
  \div_x\vecko&=0,\label{eulereq1}\\
  \rho_{\rm s}\left(\pder{\vecko}{t}+\vecko\cdot\nabla_x\vecko\right)&=\div_x\TT,\label{euler_balance}
\end{align}
where the Cauchy stress tensor reads
\begin{equation}
\TT = -p\II + 2\mu\DD + G\BB^d
\end{equation}
and the material parameters $\mu$ and $G$ differ in the fluid or solid, i.e.

\[
\label{differing_material_parameters}
    \mu=\
    \begin{cases}
    0           &\textrm{in}\ \ \Omega_{\rm s}\\
    \mu_{\rm f} &\textrm{in}\ \ \Omega_{\rm f}
  \end{cases},\qquad
  G=
  \begin{cases}
    G_{\rm s}   &\textrm{in}\ \ \Omega_{\rm s}\\
    0           &\textrm{in}\ \ \Omega_{\rm f}
  \end{cases},\qquad
  \rho=
  \begin{cases}
    \rho_{\rm s}&\textrm{in}\ \ \Omega_{\rm s}\\
    \rho_{\rm f}&\textrm{in}\ \ \Omega_{\rm f}
\end{cases}.
\]
Finally, $\BB$ satisfies
\begin{equation}
  \BB=\II\ \textrm{in}\ \Omega_{\rm f}\qquad \textrm{and}\qquad
  \pder{\BB}{t}+\vecko\cdot\nabla\BB-\LL\BB-\BB\LL^{\rm T}=\mathbb{O}\ \textrm{in}\ \Omega_{\rm s}.\label{eulereq2}
\end{equation}
In the purely Eulerian framework, displacement $\U$ is unnecessary as the whole problem is formulated only in velocity $\V$. The kinematic and dynamic interface boundary conditions are automatically satisfied for continuous velocity $\vecko$.

\mycomment{
Below, you can see the basic classification of the fluid-structure methods.
\begin{itemize}
  \item mesh-free methods
        \begin{itemize}
          \item[ex.] Smooth Particle Method (SPH)
          \item[ex.] Discrete Element Method (DEM)
        \end{itemize}
  \item mesh-based methods
        \begin{itemize}
          \item[$\bullet$] diffused boundary methods
                \begin{itemize}
                  \item[ex.] \textbf{Level-set Method} (LM)
                \end{itemize}
          \item[$\bullet$] sharp boundary methods
                \begin{itemize}
                  \item[ex.] CutCell method
                  \item[ex.] \textbf{Arbitrary Lagrangian-Eulerian} (ALE)
                \end{itemize}
        \end{itemize}
\end{itemize}
The mesh-free methods are not suitable for our problem. However, these methods could perform a rebound, but the theoretical and practical convergence are not satisfying for our purposes. Moreover, the finite number of particles could cause unnatural behavior in the area of contact.

}
\section{Arbitrary Lagrangian-Eulerian method for fluid-structure interaction}
Choosing the optimal numerical method is not a straightforward task. The choice of method depends on the character of the physical problem we would like to solve, the precision we would like to get, and the complexity of the implementation. This paper compares two approaches to the fluid-structure interaction (FSI). In this section, we employ the arbitrary Lagrangian-Eulerian (ALE) method. The following section is devoted to the Eulerian level-set approach.

The ALE method is a well-known method, first proposed by the Los Alamos group, that solved compressible Euler equations with a time-changing mesh \cite{Hirt1974}. In their method, the 
mesh configuration was always referred to the previous time step, and thus, the mesh was updated simultaneously as the simulation evolved. Let us call it the Simple ALE method. The advantage of such an approach is that it is simple to implement. However, the mesh updates after every time-step can reduce the order of accuracy. Although this can be solved by modifying this simple approach (see, for example, \cite{Hay2014} and references therein), very often, a so-called Full ALE method (sometimes called Total ALE) that does not suffer from loss of accuracy is used. In this method, the mesh configuration is kept fixed through the whole simulation, which is the main advantage of this method. Both these classical approaches were investigated/applied, for instance, in \cite{hughes1981, Donea1982, HronTurek2006a, richter2010finite, HronTurekMadlikEtAl2010, turek2011numerical, Feistaueretal2013, frei2016long}. 
For even more references, the reader is directed to two monographs on fluid-structure interaction \cite{richter2017FSIbook} and \cite{frei2018FSIbook}. 
The Full ALE method is based on constructing an arbitrary continuous unknown displacement $\U_f$ in $\Omega_{\rm f}$ (in $\Omega_{\rm s}$ the displacement $\U_s$ is the physical one) such that the function  $\chi_f = \text{id} + \U_f$ maps the computational domain into the domain deformed by the motion of the solid. In other words, we are looking for an extension of $\U_s$ in the fluid domain $\Omega_{\rm f}$. For that, we define $\U_f$ independent of the fluid-equation as a solution on the domain $\Omega_{\rm f}$, such that for every test function $\PHI$
\begin{align}
  \int_{\Omega_f(0)} \mathcal{A}(\U_f, \PHI) \, \dx &= 0, \\
  \U_f &= \U_s \; \text{at } \Gamma_{\rm int}(0), \\
  \U_f &= (0, 0) \; \text{at } \partial \Omega_{\rm f}(0) \setminus \Gamma_{\rm int}(0),
\end{align}
for some arbitrary operator $\mathcal{A}$. We choose a deformation preserving $\U_f$, such that the deformed mesh is as regular as possible. We have chosen a pseudo-solid approach \cite{SACKINGER199683}, which considers $\U_f$ as a solution to the elasticity problem with zero density. In particular, we are using an isotropic linear material

\begin{align}
  \mathcal{A}(\U_f, \PHI) = \lambda_a  \Div{\U_f}\Div{\PHI} + \mu_a (\nabla \U_f  + \nabla \U_f^T )\cdot \nabla \PHI.
\end{align}
The $\lambda_a $ and $\mu_a$ stand for Lamé coefficients, which can be suitably determined from mesh properties. We have chosen the coefficients in the following way
\begin{align}
  E_a(c) &= \frac{10}{V(c)^{\frac{9}{8}}},   &&\nu_a = - 0.02, \\
  \mu_a &= \frac{E_a}{2(1 + \nu_a)},  &&\lambda_a = \frac{\nu_a E_a}{(1 - \nu_a) (1 - 2\nu_a) }.
\end{align}
The coefficients depend on the volume $V(c)$ of a cell $c$.

\subsection{Re-meshing} \label{sec:remeshing}
We obtain the current mesh by applying the displacement to the mesh in the reference configuration. The Full ALE method works well if the deformation is not too large, which is the case of the well-known FSI benchmark \cite{HronTurek2006,HronTurek2006b}.

However, if the displacement is large enough, the regularity and topology of the current mesh can be violated. Although the Full ALE method evaluates all integrals on the mesh in the fixed configuration, the properties of the current mesh could cause the system to be ill-conditioned.
The re-meshing process is designed to deal with this problem. Here, we take advantage of both classical ALE approaches (Simple ALE and Full ALE) and employ the Updated ALE method, where the computational mesh is created by deforming the previous mesh by the displacements $u_s$ and $u_f$, respectively. The difference with respect to the Simple ALE is that the remeshing is not applied at every time step. A detailed description is given in the following subsection.

In particular, the vertices $x_i$ are shifted by 
\begin{equation}
  x_i^n = \begin{cases}
            x_i^{n-1} + u^n_f(x_i^{n-1}, t^n) \text{ in } \Omega_{\rm f}, \\
            x_i^{n-1} + u^n_s(x_i^{n-1}, t^n) \text{ in } \Omega_{\rm s}.
          \end{cases}
\end{equation}
This mesh does not need to be regular anymore. Thus, we apply an adaptive strategy afterward and 
further computations are performed on the repaired mesh.

How to create the new mesh in the re-meshing process is not a straightforward task. We would like to satisfy the following mesh properties: (i) Keep subdomains; (ii) Keep the interface edges; (iii) Small mesh changes; (iv) Compatibility with FEniCS. We have created a library ADmesh \cite{ADmesh} that satisfies all the above criteria.

One strategy to create mesh containing subdomains is based on the level-set method. In this strategy, we start with a regular mesh and move the vertices close to the interface given by the level-set function so that the interface is mimicked with edges. This strategy has been studied in \cite{cbea7e2e-fd5a-36dc-b3ae-94aa290a1a12, xie2011finite} and in \cite{fang2013isoparametric} it is done even with curved edges. This strategy would have some advantages in our case. The mesh would be regular almost everywhere in the whole domain.
Moreover, the refinement would not be complicated because it can be done only on the initial regular mesh. However, it also has some drawbacks. The vertices on the interface do not have to be in the same position after the re-meshing. Thus, some part of the fluid would be assigned into the solid part and vice versa. Further, the level-set function would need to be advected with a velocity field, which can be problematic.

Recreating the mesh from the edges of the interface and the boundary is also possible. This could be done by some existing library (Gmsh  \cite{gmsh}, Triangle \cite{Triangle} (using Ruppert’s Delaunay refinement algorithm \cite{shewchuk2002delaunay})). The newly created mesh could have even better properties than those obtained by ADmesh because of the superior algorithms in the libraries. However, the mesh could be completely different. Thus, the interpolation error would be more significant. Further, coupling the existing software with the FEniCS can be technically difficult.

Local mesh operations are a well-known technique described below for improving mesh quality, see \cite{dolejvsi1998anisotropic}, \cite{li20053d}, \cite{de1999parallel}. This has also been used in the FSI method, see \cite{compere2010mesh}. The ADmesh uses local operations such that the regularity of the mesh is iteratively improved (as described in the following subsection). The ADmesh is designed to work smoothly with FEniCS and even refine the mesh based on the functions in the FEniCS format. Further, we have been able to obtain meshes by this algorithm regularly enough, so further algorithm improvement does not seem necessary. The ADmesh allows us to satisfy all the requirements on the mesh that we have postulated before.

\begin{figure}[H]
\begin{subfigure}{.5\textwidth}
  \centering
  \includegraphics[scale=0.4]{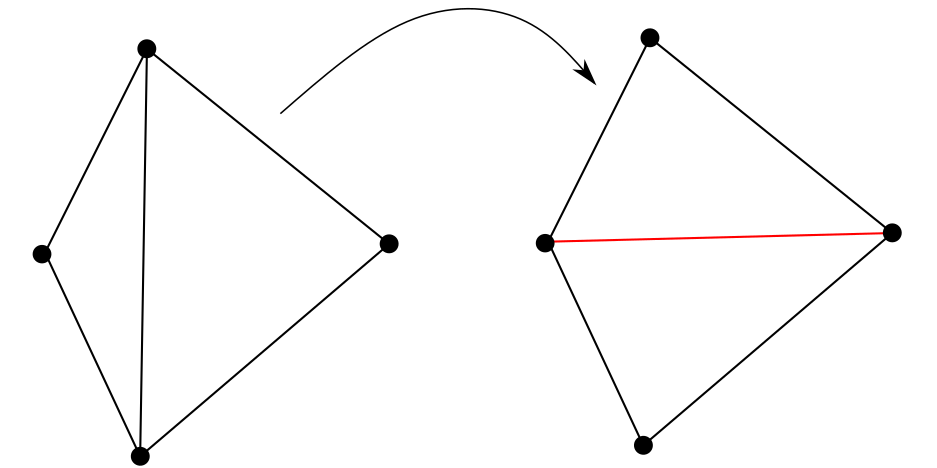}
  \caption{edge flipping}
  \label{fig:sfig1}
\end{subfigure}%
\begin{subfigure}{.5\textwidth}
  \centering
  \includegraphics[scale=0.4]{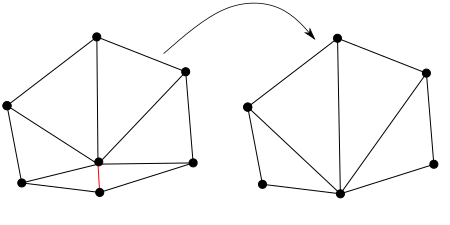}
  \caption{edge reduction}
  \label{fig:sfig2}
\end{subfigure}\\[10pt]
\begin{subfigure}{.5\textwidth}
  \centering
  \includegraphics[scale=0.4]{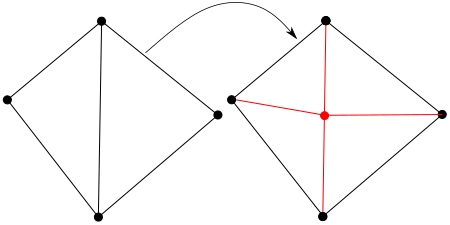}
  \caption{edge splitting}
  \label{fig:sfig3}
\end{subfigure}
\begin{subfigure}{.5\textwidth}
  \centering
  \includegraphics[scale=0.4]{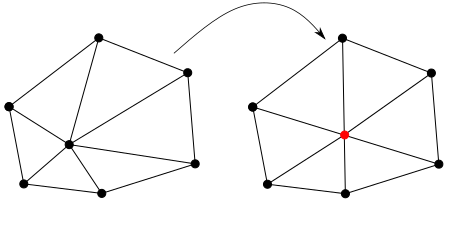}
  \caption{vertex movement}
  \label{fig:sfig4}
\end{subfigure}
\caption{Illustration of the mesh operations implemented in the ADmesh library.}
\label{fig:mesh_operations}
\end{figure}

\begin{figure}[H]
  \centering
  \includegraphics[scale=0.17]{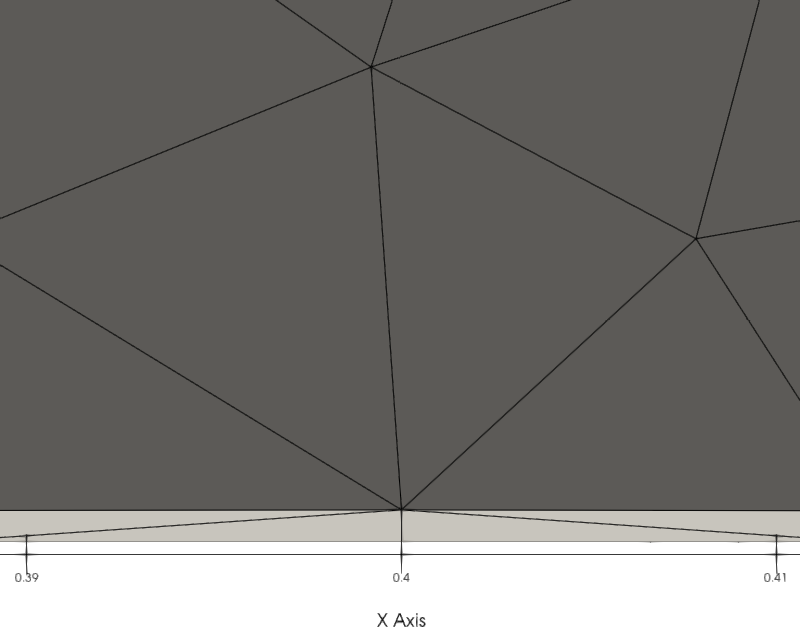}
  \includegraphics[scale=0.17]{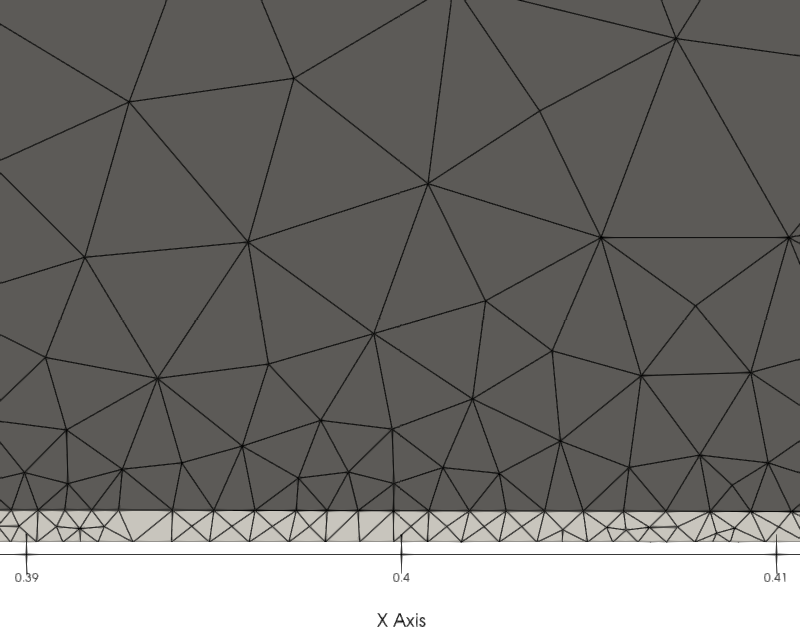}
  \includegraphics[scale=0.17]{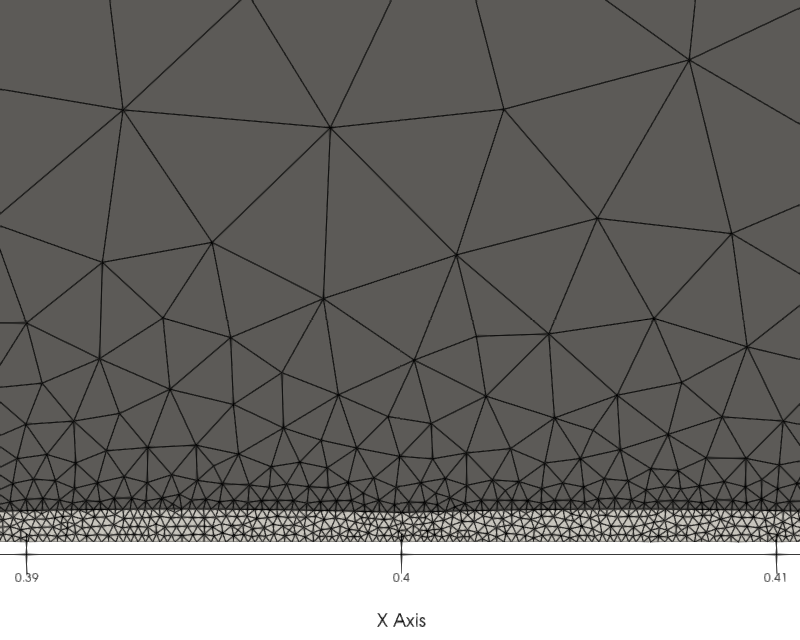}
  \caption{The figures show the difference between the refinement strategies during the rebound: $\noref$ Refinement strategy (left), $\qualityref$ Refinement strategy (middle), and $\eicref$ Refinement strategy (right).}
  \label{fig:ref_streategy_screenshot}
\end{figure}

\paragraph{Mesh adaptation}
We consider a triangular mesh with an undisrupted topology. We aim to increase the regularity of the triangles and not disrupt the edges along the interface. We have chosen the adaptive mesh technique based on several local mesh operations, such as flipping, reduction, splitting an edge, and vertex movement (see Figure~\ref{fig:mesh_operations}), which are iteratively applied wherever it increases the regularity, see \cite{dolejvsi1998anisotropic}.
The mesh adaptation can be equipped with a refinement strategy. In our case, we would like to refine the rebound area. We have tested three possible strategies described below; see Figure~\ref{fig:ref_streategy_screenshot}.
\paragraph{$\noref$ Refinement}
This method tries to keep the length of the edges during the rebound. This leads to loose of the regularity of the triangles. Eventually, this could cause a loss of convergence in Newton solver. However, the number of triangles is not increasing, and we see that the result can still be meaningful.
\paragraph{$\qualityref$ Refinement}
To keep the regularity of cells, we can refine the mesh where the quality is lower than a chosen constant. This method can ensure only one layer of cells under the ball during the rebound. The diameter of these cells will be proportional to the distance under the ball.
\paragraph{$\eicref$ Refinement}\label{refinement:Eiconal}
If we want to keep more than one layer of cells under the ball during the rebound, we must measure the distance between the solid and the wall. For every point $x\in\Omega_{\rm f}$, we define the distance to the ball $e_b$ and the distance to the wall $e_w$. The distance function $e_b$ is obtained as a solution to the $\eicref$ equation
\begin{align}
  \norm{\nabla e_b} &= 1 \text{ in } \Omega(t_r),\\
  e_b &= 0 \text{ at } \Gamma_{\rm int},
\end{align}
where $t_r$ denotes the time of re-meshing.
Then, the function $e$ computed as $e = e_w + e_b$ indicates how much refinement is needed.

\subsection{Updated ALE method}
We will study a problem where a large deformation is expected. Thus, at some time $t_r$, the mesh ceases to be regular enough, and we change the given computational domain to the current Eulerian grid $\Omega(t_r)$. The Updated ALE method has been known for some time \cite{stein2003,bungartz2006}. We describe our implementation that employs our own meshing library, ADmesh.

\begin{figure}[H]
  \centering
  \includegraphics[scale=0.4]{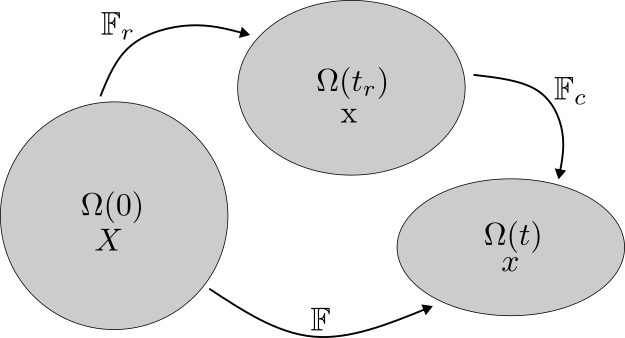}
  \caption{Multiplicative split of the total deformation gradient $\FF$ into the re-meshed part $\FF_r$ and the current part $\FF_c$.}
  \label{fig:remesh_config}
\end{figure}

After a possible re-meshing, we save the deformation gradient $\FF_r$, which describes the transformation between the reference domain $\Omega(0)$ and the computational domain $\Omega(t_r)$. Since we seek a deformation and a velocity as global functions $(\V^n, \U_c^n)$ on the whole domain $\Omega$, all calculations must be performed in the new configuration. The current domain $\Omega(t)$ is related to the computational domain $\Omega(t_r)$ by the current part of the deformation gradient $\FF_c$, see Figure~\ref{fig:remesh_config}. Thus, the total deformation gradient $\FF$ is multiplicatively split by
 \begin{equation}
\FF=\FF_c\FF_r.
\end{equation}

The current part of the deformation gradient $\FF_c$ can be expressed by introducing a new displacement $\U_c$ that measures the deformation from the computational domain $\Omega(t_r)$ to the current domain $\Omega(t)$, i.e.
\begin{equation}
    \FF_c(\mathrm{x}, t) = \II + \nabla_{\mathrm{x}} \U_c(\mathrm{x}, t).
\end{equation}

Since all equations have to be expressed in the computational domain $\Omega(t_r)$ and its boundary part $\Gamma_{\rm int}(t_r)$, the use of integral substitution theorem rewrites the volume and surface integrals:
\begin{align}
  \int_{\Omega(t)} f(x) \,\dx &= \int_{\Omega(t_r)} \det (\FF_c) f(\mathrm{x}) \,\mathrm{d}\mathrm{x} \\
  \int_{\Omega(0)} f(X) \,\mathrm{d}X &= \int_{\Omega(t_r)} \det (\FF_r^{-1}) f(\mathrm{x}) \,\mathrm{d}\mathrm{x}\\
  \int_{\Gamma_{\rm int}(t)} \F(x) \cdot \N(x, t) \,\mathrm{d}s &= \int_{\Gamma_{\rm int}(t_r)} \det (\FF_c) \F(\rx) \cdot \FF_c^{-T} \N(\rx, t_r) \,\mathrm{d}\mathrm{s} \\
  \int_{\Gamma_{\rm int}(0)} \F(X) \cdot \N(X, 0) \,\mathrm{d}S &= \int_{\Gamma_{\rm int}(t_r)} \det (\FF_r^{-1}) \F(\rx) \cdot \FF_r^{T} \N(\rx, t_r) \,\mathrm{d}\mathrm{s}\\
  \int_{\Omega(t)} \TT(x) \cdot \nabla_{x} \PHI(x) \,\dx &= \int_{\Omega(t_r)} \det(\FF_c) \TT(\rx) \FF_c^{-T} \cdot \nabla_{\rx} \PHI(\rx) \,\mathrm{d}\rx \\
  \int_{\Omega(0)} \PP(X) \cdot \nabla_{X} \PHI(X) \,\mathrm{d}X &= \int_{\Omega(t_r)} \det(\FF_r^{-1}) \PP(\rx) \FF_r^{T} \cdot \nabla_{\rx} \PHI(\rx) \,\mathrm{d}\rx
\end{align}
The formal weak formulation for the solid can be formulated as follows.
\paragraph{Solid part in $\Omega_s$}
Let $\FF_r$ be given, we are looking for functions $(\V_s, \U_{s,c})$ such that for every suitable $(\PHI_v, \PHI_u)$ it holds
\begin{equation}
  \begin{split}
    &\int_{\Omega_s(t_r)} \det (\FF_r^{-1}) \frac{\partial \V_s}{\partial t} \cdot \PHI_v \,\mathrm{d}\rx + \int_{\Omega_s(t_r)} \det(\FF_r^{-1})\PP \FF_r^T \cdot \nabla \PHI_v \, \mathrm{d}\rx \\
    & + \int_{\Gamma_{\rm int}(t_r)} \det (\FF_c) \TT_f \FF_c^{-T} \N(\rx, t_r) \cdot \PHI_v \,\mathrm{d}\rs
      + \int_{\Omega_s(t_r)} \det(\FF_r^{-1}) \left(\V_s - \frac{\partial \U_{s,c}}{\partial t}\right) \cdot \PHI_u \,\mathrm{d}\rx = 0.
  \end{split}
\end{equation}
Here, the symbol $\PP$ stands for the first Piola-Kirchhoff stress tensor, which, in the case of incompressible neo-Hookean solid, is in the form \eqref{firstPiola2}. However, this leads to a problem of pressure discontinuity across the interface $\Gamma_{\rm int}$. 
To avoid this issue, in all ALE simulations, instead of assuming an incompressible solid, we approximate the ball by a compressible neo-Hookean solid whose first Piola-Kirchhoff stress is given by
\begin{equation}
  \label{eq:comp_Piola}
  \PP = G (\FF - \FF^{-T}) + \kappa (J - 1) J \FF^{-T},
\end{equation}
where the values of elastic bulk modulus $\kappa$ is 20 times larger than the elastic shear modulus $G$ (for $G=50$ kPa, $\kappa=1$ MPa). Such choice of constitutive relation and material parameters, on the one hand, makes the pressure continuous and, on the other hand, ensures that the ball is nearly incompressible. We have verified that, in this case, the Jacobian $J$ differs from one by less than $2\times10^{-3}$.

\paragraph{Fluid part in $\Omega_f$}
In the fluid domain we require for solution $(\V_f, \U_{f,c}, p)$, that for every test function $(\PHI_v, \PHI_u, \varphi_p), \PHI_u = \PHI_v = {\bf 0} $ on $\Gamma_{\mathrm{bot}}$ it holds
\begin{equation}
  \begin{split}
    &\int_{\Omega_f(t_r)} \rho_{f} \det(\FF_c) \left( \frac{\partial\V_f}{\partial t} + (\nabla\V_f) \FF_c^{-1} \left( \V_f - \frac{\partial{\U_{f,c}}}{\partial t} \right) \right) \cdot \PHI_v \, \mathrm{d} \rx \\
    & + \int_{\Omega_f(t_r)} \det(\FF_c) \TT_f \FF_c^{-T} \cdot \nabla \PHI_v \, \mathrm{d} \rx + \int_{\Omega_f(t_r)} \mathcal{A}(\U_{f,c}, \PHI_u) \, \mathrm{d} \rx  \\
    & + \int_{\Omega_f(t_r)} \det(\FF_c) \FF_c^{-1} \V_f \cdot \nabla \varphi_p \, \mathrm{d} \rx = 0.
  \end{split}
\end{equation}
Finally, on the interface $\Gamma_{\rm int}(t_r)$, the kinematic compatibility conditions are satisfied, which state that the displacement and velocity are continuous across the interface, i.e.
\begin{equation}
  \begin{split}
    \U_{f,c} &= \U_{s,c} \text{ on } \Gamma_{\rm int}(t_r),\\
    \V_f  &= \V_s \ \ \text{ on } \Gamma_{\rm int}(t_r).
  \end{split}
\end{equation}
In fact, these conditions are satisfied automatically because the solution is sought in a continuous finite element space in the whole domain $\Omega$, see~\eqref{eq:combined} below.

\paragraph{Algorithm}
The implementation of the re-meshing process contains three main parts. First, we shift the mesh by the displacement $\U_c$, keeping the nodal values untouched. In the second part, we create a new mesh with regular triangles, as is described in \ref{sec:remeshing}.
Finally, in the last part, we interpolate all unknowns $\U_c, \V$, and $p$ to the new mesh.
We use a function {\tt create\_transfer\_matrix} from the DOLFIN   library for that.

The pseudo algorithm for re-meshing is demonstrated below.
\\
\begin{algorithm*}[H]
\SetAlgoLined
\KwData{$v^n, u_c^n, \mesh, \FF_r^n, q_0$}
$mesh_{m} = \text{move}\_\text{mesh}(\mesh, u_c^n)$\Comment*[r]{Shift the mesh by displacement $u^n$}
\If{ $\text{quality}(\mesh_{m}) <q_0$ }{
  $\FF_r^n = \FF_c^n \FF_r^n$ \Comment*[r]{Update the deformation gradient}
  $\text{change}\_{\mesh}(v^n, \mesh_m)$ \Comment*[r]{Define functions on the $\mesh_m$,}
  $\text{change}\_{\mesh}(\FF^n_r, \mesh_m)$\Comment*[r]{nodal values are the same}
  $\mesh = \text{adapt}\_\text{mesh}(\mesh_m)$\;
  $u_c^n = 0$\;
  $\FF_c^n = \II + \nabla u_c^n$\;
  $v^n = \text{interpolate}(v^n, \mesh)$ \Comment*[r]{Interpolation on repaired mesh}
  $\FF^n_r = \text{interpolate}(\FF^n_r, \mesh)$\;
}
\Return $v^n, u_c^n, \mesh, \FF^n_r, \FF_c^n$
\caption*{Re-meshing}
\end{algorithm*}
The quality of the mesh can be evaluated differently. In this paper, we consider it in the form
\begin{equation}
    \min_{\tau \in \mathcal{T}} \frac{\text{inc}(\tau)}{\text{exc}(\tau)},
\end{equation}
where $\text{inc}(\tau)$ and $\text{exc}(\tau)$ are incircle and excircle volume of a triangle $\tau$ and $\mathcal{T}$ denotes triangulation of the mesh. The interpolations of the functions are the most time-consuming part of the re-meshing process.

\subsection{Space and time discretization}
To solve the problem numerically, we discretize the equations in time and space.
\paragraph{Space discretization} We base our space approximation on the (formal) weak formulation. We use the finite element method to approximate the function spaces using a triangular mesh of the domain $\Omega$. In particular, we use linear and quadratic Lagrange finite elements.

Since we assume continuity of the displacement and velocity across the interface $\Gamma_{\text{int}}$, we use continuous finite elements throughout the domain $\Omega$. We have chosen piecewise quadratic continuous finite elements for velocity and displacement. Further, we approximate the pressure with piecewise linear continuous elements for the pressure.

The formal weak formulation can be formulated as follows. Let $\FF_r^n$ and $t_r$ be given and let $\{\tau_i\}$ denotes the triangulation of $\Omega(t_r)$, we are looking for functions $(\V^n, \U_c^n, p^n)$ such that it holds
\begin{equation}
  \label{eq:combined}
  \begin{split}
    &\int_{\Omega_s(t_r)} \rho_s \det (\FF_r^n)^{-1} \frac{\partial\V^n}{\partial t} \cdot \PHI_v^n \,\mathrm{d}\rx + \int_{\Omega_s(t_r)} \det(\FF_r^n)^{-1}\PP^n (\FF_r^n)^T \cdot \nabla \PHI_v^n \, \mathrm{d}\rx \\
    & + \int_{\Omega_s(t_r)} \det(\FF_r^n)^{-1} \left(\V^n - \frac{\partial\U^n}{\partial t} \right) \cdot \PHI_{u}^n \,\mathrm{d}\rx \\
    & + \int_{\Omega_f(t_r)} \rho_{f} \det(\FF_c^n) \left( \frac{\partial\V^n}{\partial t} + (\nabla\V^n )\FF_c^{-1} \left(\V^n - \frac{\partial\U^n}{\partial t} \right) \right) \cdot \PHI_v^n \, \mathrm{d}\rx \\
    & + \int_{\Omega_f(t_r)} \det(\FF_c)\TT_f^n (\FF_c^n)^{-T}\cdot \nabla \PHI_v^n \,\mathrm{d}\rx + \int_{\Omega_f(t_r)} \mathcal{A}(\U_f^n, \PHI_{u}^n) \, \mathrm{d}\rx  \\
    & + \int_{\Omega(t_r)}\det(\FF_c^n) (\FF_c^n)^{-1} \V_f^n \cdot \nabla \varphi_p^n \,\mathrm{d}\rx = 0,
  \end{split}
\end{equation}
for every $(\PHI_v^n, \PHI_u^n, \varphi_p^n)$, where the functions and test functions belong to the following finite element spaces
\begin{align}
  \V^n, \PHI_v^n &\in \{ f \in P^2(\tau_i);\; f \in C(\Omega); \; f = (0, 0) \text{ at }  \Gamma_{\text{bot}}\}, \\
  \U_c^n &\in \{ f \in P^2(\tau_i);\; f \in C(\Omega); \; f = (0, 0) \text{ at }  \partial \Omega \}, \\
  \PHI_{u}^n &\in \{ f \in P^2(\tau_i);\; f \in C(\Omega); \tau_i \subset \Omega_s; \; f = (0, 0) \text{ at }  \Gamma_{\text{bot}} \cup \Gamma_{\text{int}} \}, \\
  p^n, \varphi_p^n &\in \{ f \in P^1(\tau_i);\; f \in C(\Omega); \tau_i \subset \Omega_f \},
\end{align}

\paragraph{Time discretization}\label{Sec:Glowinski}
We choose $N$ timesteps $t_{i}$ uniformly distributed in the interval $(0, T)$ where we evaluate the solution.
The most common time discretization is the backward Euler scheme, which reads
\begin{equation}\label{backward_Euler}
    u^{n + 1}(x) - u^{n}(x) = \Delta t f(u^{n + 1}).
\end{equation}
Since this scheme is only of the first order, we use the Glowinski scheme \cite{glowinski2003} that consists of two implicit steps and one explicit step, i.e.
\begin{equation}\label{eq:time_disc}
  \begin{split}
    u^{n + \theta}(x) - u^{n}(x) &= \theta \Delta t f(u^{n + \theta}),\\
    u^{n + 1 - \theta} &= \frac{1 - \theta}{\theta} u^{n + \theta} + \frac{2\theta - 1}{\theta} u^{n},\\
    u^{n + 1}(x) - u^{n + 1 - \theta}(x) &= \theta \Delta t f(u^{n + 1}),
  \end{split}
\end{equation}
where $\theta=1/\sqrt{2}$. The pressure $p$ is taken at the current time level only (similarly as in the case of the Crank-Nicolson time scheme), and no additional pressure postprocessing step (as in \cite{Turek2006}) is applied. Glowinski \cite{glowinski2003} proved that this scheme is second-order accurate; however, in the numerical experiments (\cite{hrratu2014}), the scheme exhibits even a third-order experimental order of convergence.

\subsection{Domain Discretization}
The discretization of the ball is different in ALE and the level-set method. The level-set function can be designed to describe the subdomains precisely. In the case of the ALE method, we are using straight edges to approximate the interface, so the subdomains do not precisely describe the shape of the subdomains. Although some methods exist for dealing with this issue, we do not include them in our simulations. The domain discretization, however, plays an essential role in convergence. Thus, we specify the number of lines approximating the interface.

%%%
\section{Simulation results -- numerics}
In this section, we present the simulation of the following FSI problem using the ALE approach. This section is devoted to the numerical aspects of the simulation such as the convergence with respect to time and space, as well as the required CPU time. The following section is then devoted to the physical aspects of the simulation.

\subsection{FSI problem description} \label{sec41}

We consider a FSI problem that captures an elastic ball described by a neo-Hookean solid (the constitutive relation for the first Piola-Kirchhoff stress is given by \eqref{eq:comp} for a compressible or \eqref{firstPiola2} for an incompressible solid) submerged in an incompressible Newtonian fluid (the constitutive relation for the Cauchy stress is given by \eqref{CauchyNewtonian}). At $t=0$ the ball is thrown in a gravity-free field against a bottom wall with an initial velocity of $0.5$\,m/s. 

\begin{figure}[h]
\begin{center}
\includegraphics[width=5cm]{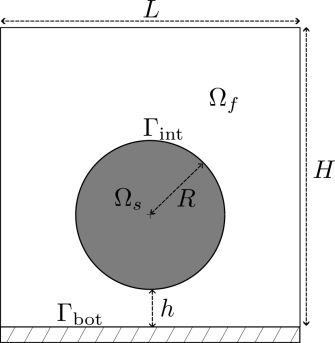}
  \caption{Problem description.}
  \label{fig:geometry}
\end{center}
\end{figure}

The problem is described in the whole domain $\Omega=(0, 0.8)\times(0, 0.8)$\,m$^2$ that is split into $\Omega_{\rm s}$ and $\Omega_{\rm f}$, see Figure~\ref{fig:geometry}. The elastic ball initially occupies a circular domain $\Omega_{\rm s}$ whose center is located in $(0.4, 0.3)$\,m with radius $0.2$\,m, so the initial distance of the ball from the wall is equal to 0.1\,m. The Newtonian fluid occupies domain $\Omega_{\rm f}$.

We assume a no-slip boundary condition at the bottom wall $\Gamma_{\rm bot}$. To capture a rebound in a much larger domain than the real computational domain, zero traction is prescribed on the remaining walls (the lateral and the top). On the interface $\Gamma_{\rm int}$, the classical kinematic (the continuity of velocities), as well as the classical dynamic (the continuity of tractions) interface conditions, are prescribed.

We assume that the fluid density is small $\rho_{\rm f}=1$~kg/m$^3$, and the solid density equals 1000~kg/m$^3$. The problem is solved for three values of viscosities 0.1~Pa\,s, 0.01~Pa\,s, and 0.001~Pa\,s (the last one corresponds to the viscosity of water). The elastic shear moduli of the solid are equal to 5, 50, and 500 kPa.

The smallness of the density has just practical reasons. The respective low fluid inertia leads to a
reasonably short contact time, which makes the problem well-computable.

\subsection{Quantities of interest}
We provide the definitions of several quantities that we record during the simulation. First six quantities minimum of the ball surface $y_{\rm min}$, minimum of the ball surface at the centre $y_{\rm min, c}$, pressure at center on the bottom $p_{\rm bc}$, kinetic energy of the ball $E_{\rm k, s}$, elastic energy of the ball $E_{\rm el, s}$ and total energy of the ball $E_{\rm s}$ are functions of time $t$, and they are defined as
\begin{align}
  y_{\rm min, c}&=\min_{(0.4,y)\in \Omega_{\rm s}}y, &&&
                                                         y_{\rm min}&=\min_{(x,y)\in \Omega_{\rm s}}y,\\
  p_{\rm bc}&=p([0.4, 0.0], t), &&&
                                    E_{\rm k, s}&=\int_{\Omega_{\rm s}} \frac{\rho_{\rm s}}{2}|\V|^2\, \dx, \\
  E_{\rm el, s}&=\int_{\Omega_{\rm s}} \frac{G}{2}\left(\tr (\BB) - 2\right) \, \dx, &&&
                                                                                         E_{\rm s}&=E_{\rm k, s}+E_{\rm el, s}.
\end{align}
Last three quantities maximum of $p_{\rm bc}$ over time, maximum of $E_{\rm el, s}$ over time and minimum of $E_{\rm k, s}$ over time are real numbers depending on the given viscosity $\mu$ and elastic modulus $G$, and they are denoted by
\begin{equation}
  \max_{t} p_{\rm bc},\quad \max_{t} E_{\rm el, s},\quad \min_{t} E_{\rm k, s}.
\end{equation}

As the initial conditions for the solid, we set $\V_s^0(x) = \V_s(0, x) = (0.0, -0.5)$. 
We must satisfy the boundary conditions and velocity continuity along the interface for the fluid part. For that purpose, we find the solution $ \V_f^0 = \V(0, \cdot) $ to Stokes equations in the fluid part, i.e.
\begin{equation}
    \begin{split}
        -\Delta \V_f^0  - \nabla p(0, x) &= {\bf 0} \text{ in } \Omega_{\rm f}, \\
        \V_f^0 &= {\bf 0} \text{ at }\Gamma_{\rm bot}, \\
        (\nabla \V_f^0 - p(0, x) \II)\N &= {\bf 0} \text{ at } \partial \Omega \setminus \Gamma_{\rm bot}, \\
        \V_f^0 &= \V_s^0 \text{ at } \Gamma_{\rm int}.
    \end{split}
\end{equation}

\subsection{Convergence study}
For a fixed viscosity $\mu=0.1$~Pa\,s and fixed elastic modulus $G=50$~kPa, we perform a convergence study to verify that our simulations converge. We use several different meshes. The meshes are denoted by $\meshlev{i}{j}$. The $i \in \{0, 1, 2, 3\}$ denotes how fine the mesh is, as it is  shown in Table~\ref{tab:cells}. The triangulation can not exactly capture the shape of the ball. Thus, we use a regular polygon to approximate it. The upper index $j \in \{100, 200, 300, 400 \}$ represents the number of vertices of the polygon.

\begin{table}[H]
  \centering
  \begin{tabular}{c|c|c|c|c|c}
    mesh    & $\meshlev{0}{200}$    & $\meshlev{1}{200}$    & $\meshlev{2}{200}$    & $\meshlev{3}{200}$ & $\meshlev{3}{200}$  \\ \hline
    $\# cells$    & $7956$        & $15165$        & $25 046$       & $37 550$ & $52153$     \\
  \end{tabular}
  \caption{Number of cells in the meshes.}
  \label{tab:cells}
\end{table}

\subsubsection{Convergence with respect to time}
\begin{figure}[H]
  \centering
  \includegraphics{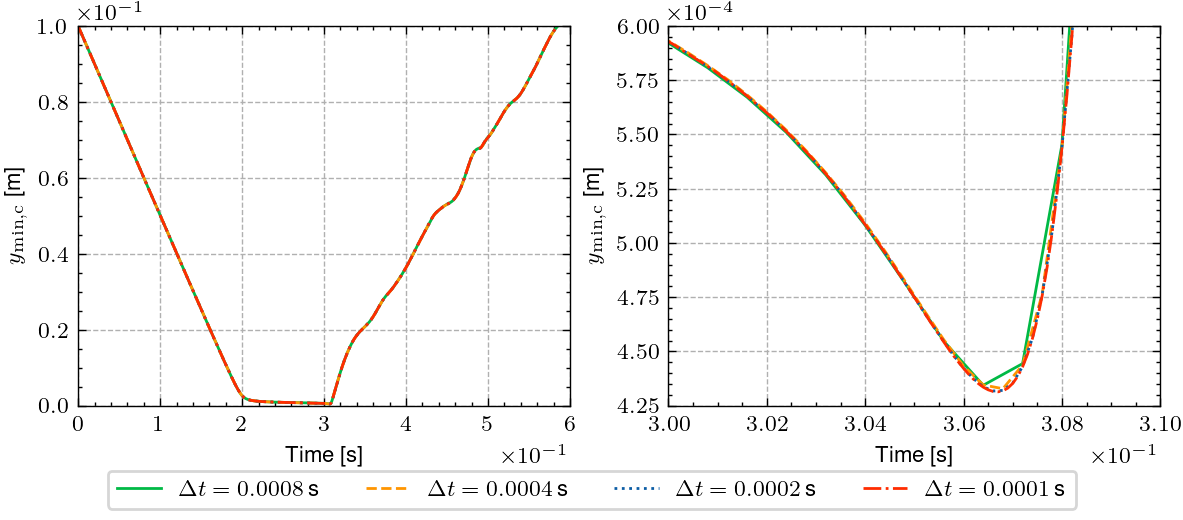}
  \caption{Comparison of time evolution of $y_{\rm min, c}$ for different timesteps $\Delta t$. The results are obtained using the viscosity $\mu=0.1$~Pa\,s, elastic modulus $G = 50$~kPa, $\meshlev{3}{600}$, and $\eicref$ strategy.}
  \label{fig:convtime}
\end{figure}

First, we check how the minimum of the ball surface $y_{\rm min}$ depends on time $t$ when decreasing the time step; see Figure \ref{fig:convtime}.
With a decreasing time step $\Delta t$, the solution converges; see Table~\ref{tab:time_conv}.

\begin{table}[!htbp]
  \centering
  \begin{tabular}{c|c|c|c|c}
    $\Delta t [s]$                     & $8\times 10^{-4}$    & $4\times 10^{-4}$    &  $2\times 10^{-4}$   & $1\times 10^{-4}$    \\ \hline
    $\text{min}_t y_{\rm min, c}$  [m]& $4.361 \times 10^{-4}$    & $4.338 \times 10^{-4}$    &  $4.332\times 10^{-4}$   & $4.330\times 10^{-4}$  \\
    $\max_{t} p_{\rm bc}$  [Pa]& $23068.022 $    & $ 23070.493$    &  $23068.551 $   & $23065.842 $  \\
    $\max_{t} E_{\rm el, s}$  [J]& $11.220$    & $11.220$    &  $11.220$   & $11.220$  \\
    $\min_{t} E_{\rm k, s}$  [J]& $8.966 \times 10^{-2}$    & $8.760\times 10^{-2}$    &  $8.757\times 10^{-2}$   & $8.755\times 10^{-2}$  \\
  \end{tabular}
  \caption{Comparison of results for different timesteps $\Delta t$. The results are obtained using the elastic modulus $G = 50$~kPa, viscosity $\mu=0.1$~Pa\,s, $\meshlev{3}{600}$, and $\eicref$ strategy.}
  \label{tab:time_conv}
\end{table}

\subsubsection{Convergence with respect to space}
\begin{figure}[H]
  \centering
  \includegraphics{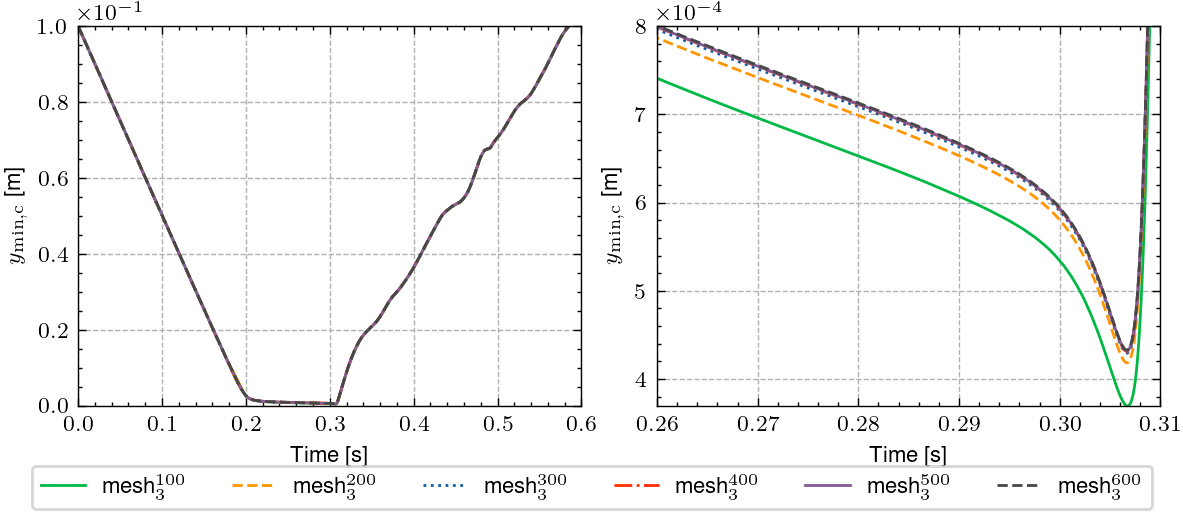}
  \caption{Comparison of time evolution of $y_{\rm min, c}$ for different discretization of the elastic ball. The results are obtained using the elastic modulus $G=50$~kPa, viscosity $\mu = 0.1$~Pa\,s, $\meshlev{3}{i}$, timestep $\Delta t = 10^{-4}$~s and $\eicref$ refinement.}
  \label{fig:points_covergence}
\end{figure}

Further, we study how the $y_{\rm min}$ evolves with the change of space resolution. We split this issue into three steps. In the first step, we study the discretization of the solid. Since we approximate the ball with a polygon, we study how the rebound behaves with a different number of vertices on the ball's surface. The results are very close for the finest discretization of the ball as it is captured in Figure~\ref{fig:points_covergence} and Table~\ref{tab:point_convergence}.

\begin{figure}[H]
  \centering
  \includegraphics{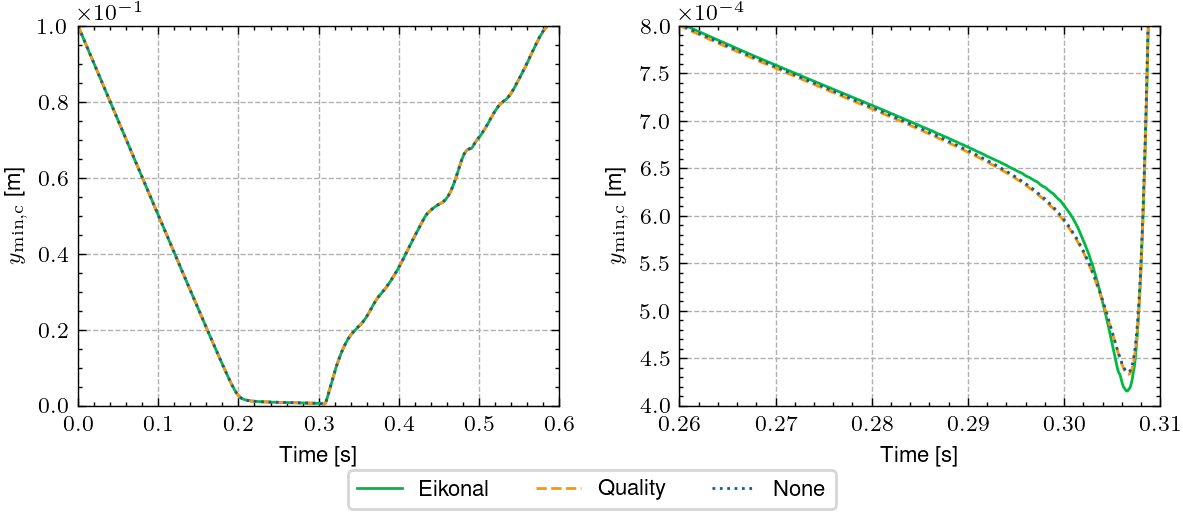}
  \caption{Comparison of time evolution of $y_{\rm min, c}$ for different refinement strategies. The results are obtained using the viscosity $\mu = 0.1$~Pa\,s, $\meshlev{3}{600}$ and timestep $\Delta t = 10^{-4}$~s.}
  \label{fig:refin_convergence}
\end{figure}

\begin{table}[H]
{\footnotesize
  \centering
  \begin{tabular}{c|c|c|c|c|c|c}
    mesh      & $\meshlev{3}{100}$    & $\meshlev{3}{200}$   & $\meshlev{3}{300}$   & $\meshlev{3}{400}$  & $\meshlev{3}{500}$   & $\meshlev{3}{600}$    \\ \hline
    $\text{min}_t y_{\rm min, c}$  [m]
    & $3.696 \times 10^{-4}$
    & $4.185\times 10^{-4}$
    & $4.294\times 10^{-4}$
    & $4.313\times 10^{-4}$
    & $4.322\times 10^{-4}$
    & $4.330\times 10^{-4}$ \\
    $\max_{t} p_{\rm bc}$  [Pa]
    & $23062.602 $
    & $23069.368 $
    & $23066.973 $
    & $23066.280 $
    & $23069.704 $
    & $23065.842 $  \\
    $\max_{t} E_{\rm el, s}$  [J]&
    $11.212$    
    & $11.218$
    & $11.218$
    & $11.219$
    & $11.220$
    & $11.220$  \\
    $\min_{t} E_{\rm k, s}$  [J]
    & $8.740 \times 10^{-2}$    
    & $8.751\times 10^{-2}$    
    & $8.747\times 10^{-2}$   
    & $8.754\times 10^{-2}$  
    & $8.755\times 10^{-2}$   
    & $8.755\times 10^{-2}$  \\
  \end{tabular}
  \caption{
    Comparison of results for different discretization of the elastic ball. The results are obtained with the elastic modulus $G = 50$~kPa, viscosity $\mu=0.1$~Pa\,s, timestep $\Delta t = 10^{-4}$~s, $\meshlev{3}{i}$, and $\eicref$ strategy.
  }
  \label{tab:point_convergence}}
\end{table}
The important question is whether the solution is unaffected by the small number of elements between the solid and the rigid wall. This issue is studied in the second step, where we focus on the refinement strategies: $\noref$ refinement, $\qualityref$ refinement, and $\eicref$ refinement strategy, see Section~\ref{refinement:Eiconal}. The $\eicref$ approach provides the most trustworthy results because the area between the wall and the ball contains significantly more elements. However, it makes the problem more expensive during the rebound regarding the number of degrees of freedom (DOFs), as shown in Table~\ref{tab:dofs}. The results obtained by the $\eicref$ strategy greatly agree with the $\qualityref$ strategy, as shown in Figure~\ref{fig:refin_convergence} and Table~\ref{tab:refinement}. The macroscopic behavior is almost the same for $\noref$ refinement. However, the $y_{\rm min, c}$ is noticeably lower during the rebound. More importantly, for lower viscosities, the $\noref$ refinement strategy fails (as the elements are getting too flat when approaching the rebound, which violates the solver's convergence). It seems quite noticeable that the Eikonal refinement strategy, allowing viscosities as low as $0.001$ Pa s, does not need significantly more elements for $\mu=0.1$ Pa s than the other methods. 

\begin{table}[H]
  \centering
  \begin{tabular}{c|c|c|c}
    strategy &  $\noref$ &  $\qualityref $ & $\eicref$ \\ \hline
    \# DOFs  at $t=0.25$ s & 299605 & 323710 & 371960 \\
    the mean value of \# DOFs & 303677 & 314507& 322425 \\
  \end{tabular}
  \caption{Number of degrees of freedom for different refinement strategies using the viscosity $\mu=0.1$ Pa s on $\meshlev{3}{600}$.}
  \label{tab:dofs}
\end{table}

\begin{table}[H]
  \centering
  \begin{tabular}{c|c|c|c}
    refinement    & $\noref$    & $\qualityref$    & $\eicref$  \\ \hline
    $\text{min}_t y_{\rm min, c}$  [m]& $4.153 \times 10^{-4}$        & $4.352 \times 10^{-4}$        & $4.330 \times 10^{-4}$    \\
    $\max_{t} p_{\rm bc}$  [Pa]& $23105.830 $    & $23069.439 $    &  $23065.842 $ \\
    $\max_{t} E_{\rm el, s}$  [J]& $11.219$    & $11.221$    &  $11.220$   \\
    $\min_{t} E_{\rm k, s}$  [J]& $8.737 \times 10^{-2}$    & $8.755\times 10^{-2}$    &  $8.755\times 10^{-2}$  \\
  \end{tabular}
  \caption{Comparison of results for different refinement strategies. The results are obtained using the viscosity $\mu=0.1$~Pa\,s, elastic modulus $G = 50$~kPa, $\meshlev{3}{600}$ and timestep $\Delta t = 10^{-4}$~s.}
  \label{tab:refinement}
\end{table}

\begin{figure}[!htbp]
  \centering
  \includegraphics{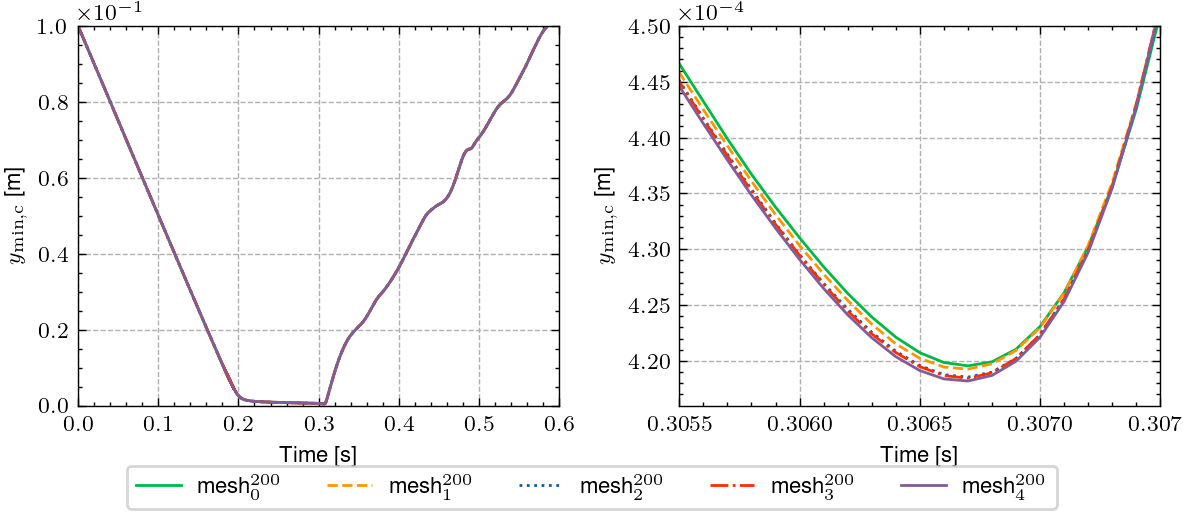}
  \caption{Comparison of time evolution of $y_{\rm min, c}$ using different mesh discretization. The results are obtained using the viscosity $\mu = 0.1$~Pa s, the elastic modulus $G = 50$~kPa, timestep $\Delta t = 10^{-4}$~s, $\meshlev{i}{200}$ and $\eicref$ refinement.}
  \label{fig:h_covergence}
\end{figure}

In the third step, we focus on convergence with respect to the global refinement. For this purpose, we fix the discretization of the ball, and we study the convergence with respect to space discretization. Since it is difficult to construct and refine the $\meshlev{0}{j}$ for large $j$ due to the large difference in the edge sizes, we use $j = 200$ for this convergence study. As shown in Figure~\ref{fig:h_covergence} and Table~\ref{tab:h_convergence}, the problem is not sensitive to the cell size.

\begin{table}[H]
  \centering
  \begin{tabular}{c|c|c|c|c|c}
    mesh
    & $\meshlev{0}{200}$
    & $\meshlev{1}{200}$
    & $\meshlev{2}{200}$
    & $\meshlev{3}{200}$
    & $\meshlev{4}{200}$
    \\ \hline
    $\text{min}_t y_{\rm min, c}$ [m]
    & $4.196 \times 10^{-4}$
    & $4.193 \times 10^{-4}$
    & $4.186 \times 10^{-4}$
    & $4.185 \times 10^{-4}$
    & $4.182 \times 10^{-4}$
    \\
    $\max_{t} p_{\rm bc}$  [Pa]
    & $23041.000$
    & $23050.815$
    & $23060.702$
    & $23069.368$
    & $23072.092$
    \\
    $\max_{t} E_{\rm el, s}$  [J]
    & $11.206$
    & $11.212$
    & $11.217$
    & $11.218$
    & $11.219$
    \\
    $\min_{t} E_{\rm k, s} $  [J]
    & $8.704  \times 10^{-2}$
    & $8.729  \times 10^{-2}$
    & $8.743  \times 10^{-2}$
    & $8.751  \times 10^{-2}$
    & $8.754  \times 10^{-2}$
    \\
 \end{tabular}
  \caption{Comparison of results for different refinement of the mesh $\meshlev{i}{200}$. The results are obtained using the viscosity $\mu=0.1$~Pa\,s, elastic modulus $G = 50$~kPa, timestep $\Delta t = 10^{-4}$~s, and $\eicref$ strategy.}
  \label{tab:h_convergence}
\end{table}
The results have been obtained on CPUs connected into a node consisting of 2x Intel(R) Xeon(R)  Gold 6140 CPU  @ 2.30GHz (36 cores in total) with 128GB RAM. The following table shows the mean time to solve one timestep. Because we use the Glowinski time scheme, each timestep consists of two nonlinear problems solved by the Newton method.

\begin{table}[H]
    \centering
    \begin{tabular}{c|c|c}
    mesh  & mean solver time [s] & mean number of DOFs \\ \hline
    $\mesh_{0}^{200}$ & $1.57$ & $89941 $ \\
    $\mesh_{1}^{200}$ & $2.28$ & $141053$ \\
    $\mesh_{2}^{200}$ & $3.55$ & $204707$ \\
    $\mesh_{3}^{200}$ & $5.33$ & $307634$ \\
    $\mesh_{4}^{200}$ & $7.24$ & $387226$ \\
    \end{tabular}
    \caption{Mean CPU times to solve one timestep of the problem for different resolutions. The results are obtained using the viscosity $\mu = 0.1$~Pa s, the elastic modulus $G = 50$~kPa, timestep $\Delta t = 10^{-4}$~s, $\meshlev{i}{200}$ and $\eicref$ refinement.}
    \label{tab:my_label}
\end{table}
%%%
The problem is computed for time $t$ upto 0.6 s. For $\Delta t=10^{-4}$ s we need to perform 6\,000 timesteps which in case of $\mesh_4^{200}$ takes approximately 12 hours of CPU time.

\section{Simulation results -- physics}
In this section, we study the simulation results from the point of view of physics.
We observe how the results change as the viscosity $\mu$ is decreased and/or how they depend on the elastic modulus $G$. In all cases, we observe a contactless rebound of the ball from the wall.

\subsection{Numerical results for different viscosities (and vacuum case) and elastic moduli}
\begin{figure}[H]
  \centering
  \includegraphics{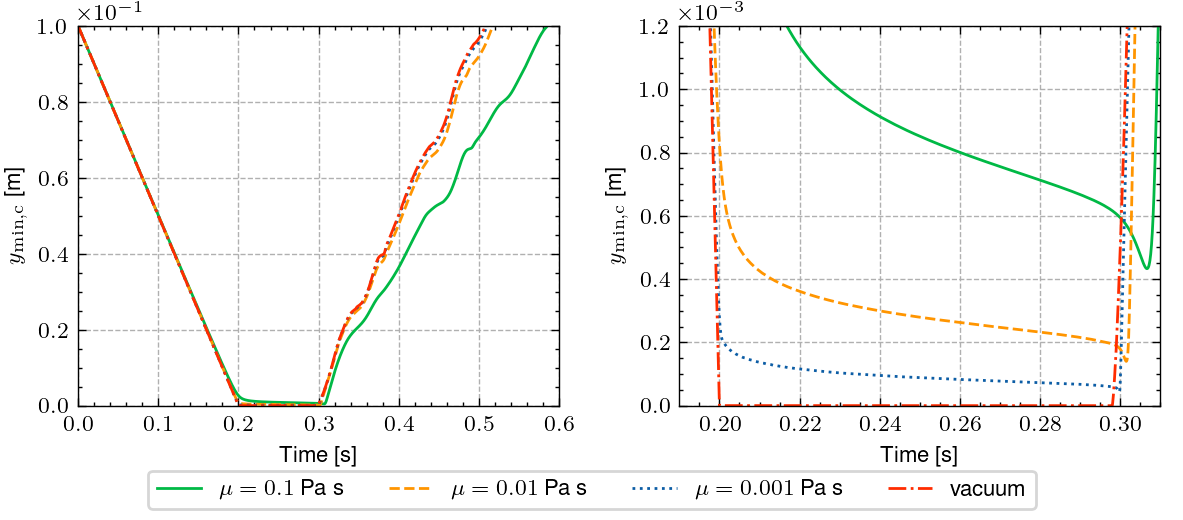}
  \caption{Comparison of time evolution of $y_{\rm min, c}$ using different viscosities. The values are computed with elastic modulus $G = 50$~kPa, timestep $\Delta t = 10^{-4}$~s, $\meshlev{3}{600}$ and $\eicref$ refinement. The rebound in a vacuum is computed with timestep $10^{-4}$~s and on a circle mesh with $600$ vertices on the boundary.}
  \label{fig:depviscosity}
\end{figure}

\begin{figure}[H]
  \centering
  \includegraphics{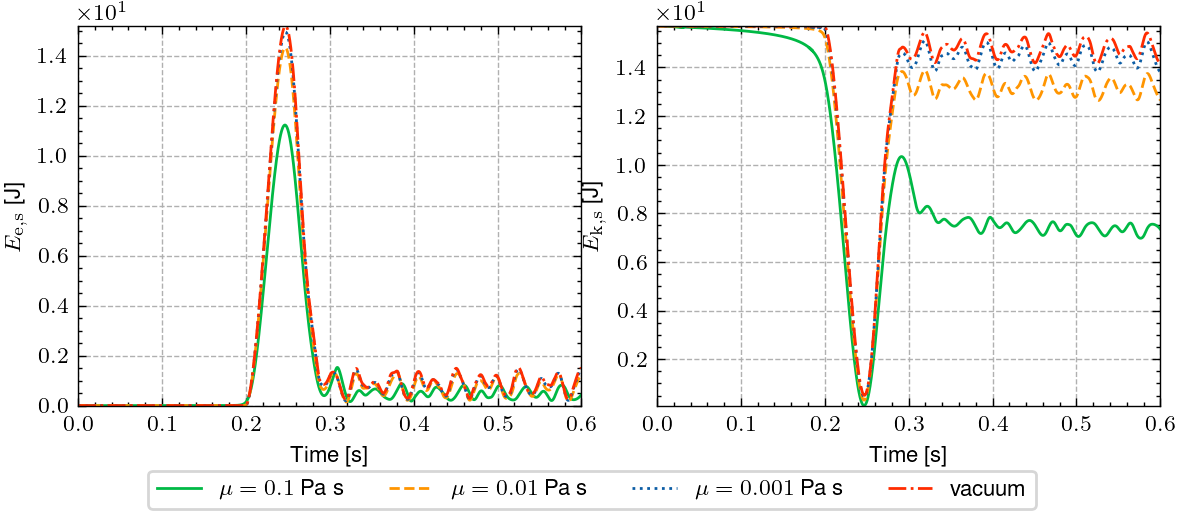}
  \caption{Comparison of time evolution of $E_{\rm e, s}$ and $E_{\rm k, s}$ using different viscosities. The results are obtained using the elastic modulus $G = 50$~kPa, timestep $\Delta t = 10^{-4}$~s, $\meshlev{3}{600}$ and $\eicref$ strategy.}
  \label{fig:depviscosity_energy}
\end{figure}

We check how the results change with the choice of viscosity $\mu$, see Table~\ref{tab:mu_convergence}, Figures~\ref{fig:depviscosity_energy} and \ref{fig:depviscosity}.
For comparison, we also compute the rebound of the elastic ball from the wall in a vacuum. In this case, the contact is implemented using the augmented Lagrangian method. The details are available in \ref{Appendix_vacuum}.

Clearly, if the ball was in the vacuum, it would reach the wall in 0.2~s (initial distance 0.1~m divided by initial velocity 0.5~m/s). With increasing viscosity, the drag increases, and the ball hits the wall later. However, it does not touch the wall but bounces off due to a pressure singularity at a finite distance from the wall. This distance is lower as the viscosity is decreased and with the lowest value $y_{\rm min}=2.71\times 10^{-5}$\,m for the lowest viscosity $\mu=0.001$~Pa\,s.
Finally, for a relatively long time around 0.1\,s (rebound phase), the bottom of the ball remains at a similar distance from the wall, and then after $t=0.3$\,s, it bounces off.

\begin{figure}[H]
  % \centering
  \includegraphics[valign=t]{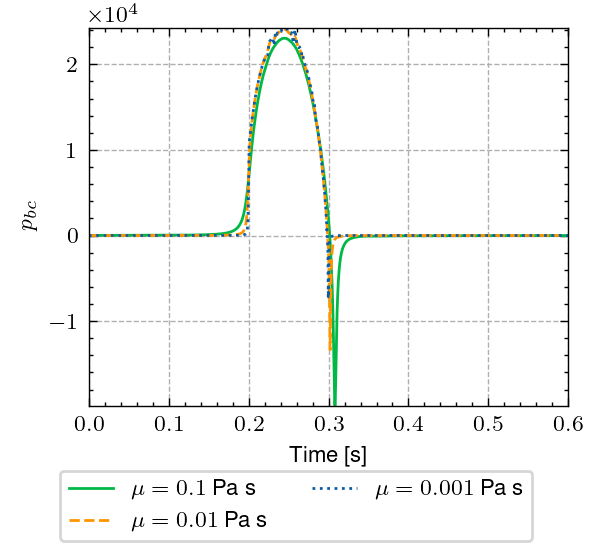}
  \includegraphics[valign=t]{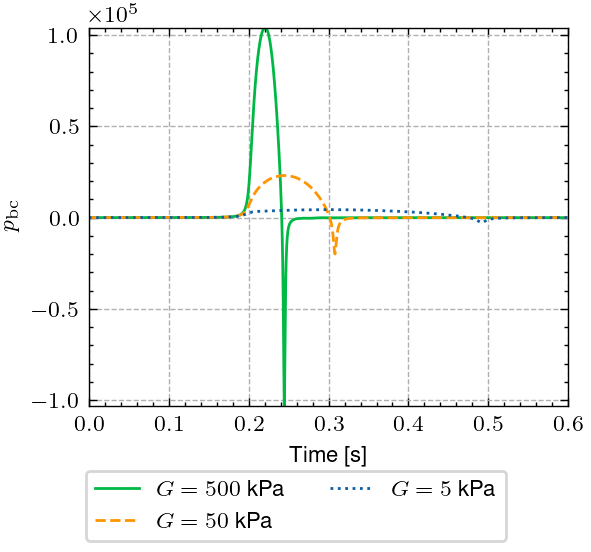}
  \caption{Evolution of the pressure $p_{\rm bc}$ with respect to time for: different viscosities $\mu = 0.1, 0.01, 0.001$~Pa\,s with fixed elastic modulus $G = 50$ kPa (left); different elastic moduli $G=5, 50, 500$~kPa with fixed viscosity $\mu=0.1$~Pa\,s (right). The results are obtained using timestep $\Delta t = 10^{-4}$~s, $\meshlev{3}{600}$ and $\eicref$ strategy.}
  \label{fig:pressure-bmcbm}
\end{figure}

The dependence of the pressure $p_{\rm bc}$ with respect to time for different viscosities are plotted in Figure~\ref{fig:pressure-bmcbm}. The pressure increases rapidly after $t=0.2$\,s when the ball reaches the wall, it reaches its maximal value at $t=0.245$\,s and then it `falls down' to the negative values at $t=0.301$\,s which enables the ball the bounce of from the wall. It is worth noticing that the maximal value of $p_{\rm bc}$ is not significantly changing with the choice of viscosity.

\begin{table}[H]
  \centering
  \begin{tabular}{c|c|c|c|c}
    $\mu$ [Pa s]    & $0.1$    & $0.01$    & $0.001$  & vacuum$^*$ \\ \hline
    $\text{min}_t y_{\rm min, c}$ [m] & $ 4.330 \times 10^{-4}$  & $1.398\times 10^{-4}$  & $4.334 \times 10^{-5}$  &  $ - $\\
    $\max_{t} p_{\rm bc}$  [Pa]& $23065.842 $    & $24113.476 $    &  $24311.393$  \\
    $\max_{t} E_{\rm el, s}$  [J]& $ 11.220 $    & $14.311 $    &  $ 14.939 $ &  $15.088$\\
    $\min_{t} E_{\rm k, s}$  [J]& $ 8.755\times 10^{-2}$    & $3.124\times 10^{-1}$   &  $4.448\times 10^{-1}$ & $5.691 \times 10^{-1}$\\
  \end{tabular}
  \caption{Comparison of the results for different viscosities $\mu = 0.1, 0.01, 0.001 $~Pa\,s and vacuum. The results are obtained using the elastic modulus $G = 50$~kPa, the timestep $\Delta t = 10^{-4}$~s, $\meshlev{3}{600}$ and $\eicref$ strategy.}
  \label{tab:mu_convergence}
\end{table}

By comparing the results, we observe that with decreasing viscosity, the elastic energy $E_{e, s}$, the kinetic energy $E_{k, s}$ and the $y_{\rm min, c}$ converge to the vacuum case.
This ensures that the computations are meaningful and strengthens the hypothesis that the non-contact rebound could be a physical phenomenon.
It is important that the density is small. Otherwise, the convergence to the rebound in the vacuum would not be satisfied.
Although the maximum pressure $p_{\rm bc}$ does not depend significantly on the viscosity, it depends greatly on the elastic modulus $G$, see Figure~\ref{fig:pressure-bmcbm}.

\begin{figure}[H]
  \centering
  \includegraphics[width=\textwidth]{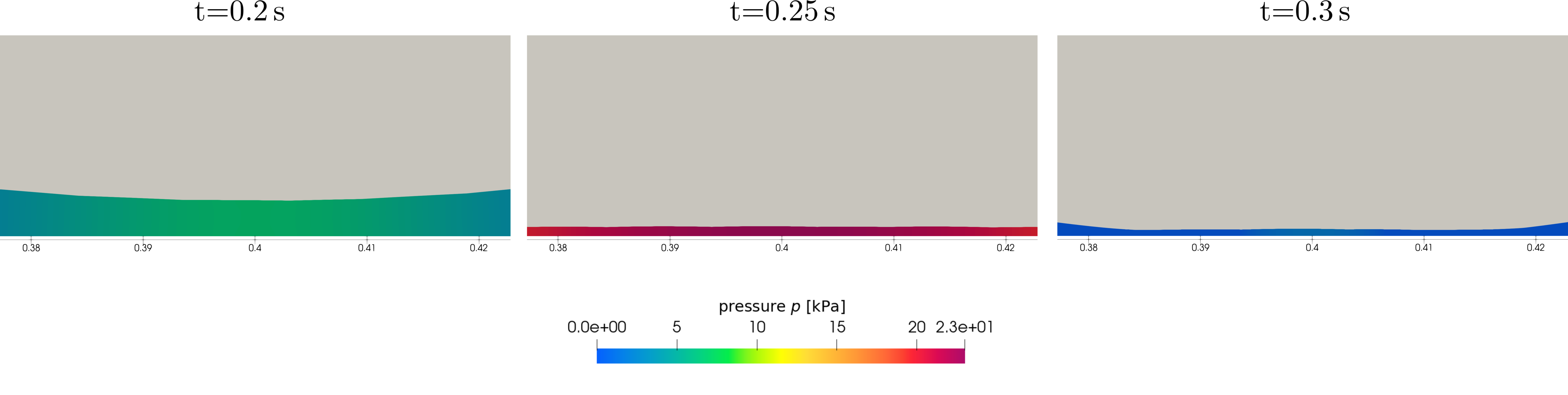}
  \caption{The pressure under the ball captured at times $t = 0.2, 0.25$ and $0.3$~s for: $\mu = 0.1$~Pa\,s, elastic modulus $G = 50$~kPa.}
  \label{fig:snap_ball_pressure}
\end{figure}

\begin{figure}[H]
  \centering
  \includegraphics{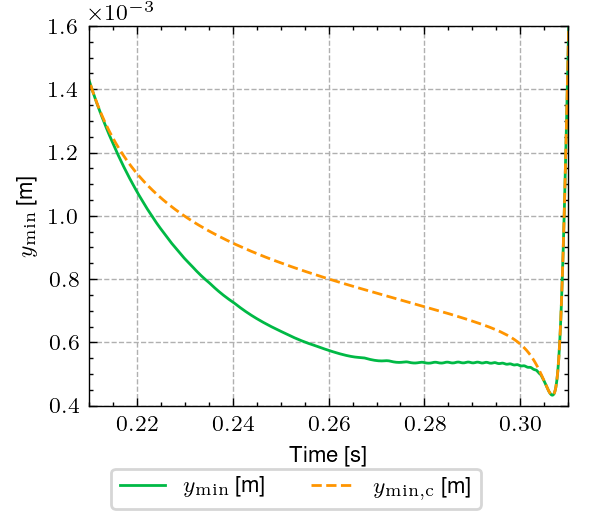}
  \caption{Comparison of time evolution of $y_{\rm min}$ and $y_{\rm min, c}$ for viscosity $\mu=0.1$~Pa\,s, elastic modulus $G=50$~kPa, timestep $\Delta t = 10^{-4}$ s, $\meshlev{3}{600}$ and $\eicref$ strategy.}
  \label{fig:mid_vs_min}
\end{figure}

\paragraph{Non-convexity of the ball during the rebound phase}
The shape of the bottom surface of the ball changes significantly during the rebound phase. See a snapshot from the simulation in Figure~\ref{fig:snap_ball_pressure}. The ball ceases to be convex, which is observed in Figure~\ref{fig:mid_vs_min} that shows the dependence of $y_{\rm min}$ and $y_{\rm min, c}$ on time $t$. Clearly, whenever $y_{\rm min}$ is lower than $y_{\rm min, c}$, the ball is not convex.

\subsection{Coefficient of restitution}

\begin{figure}[H]
  \centering
  \includegraphics[valign=t]{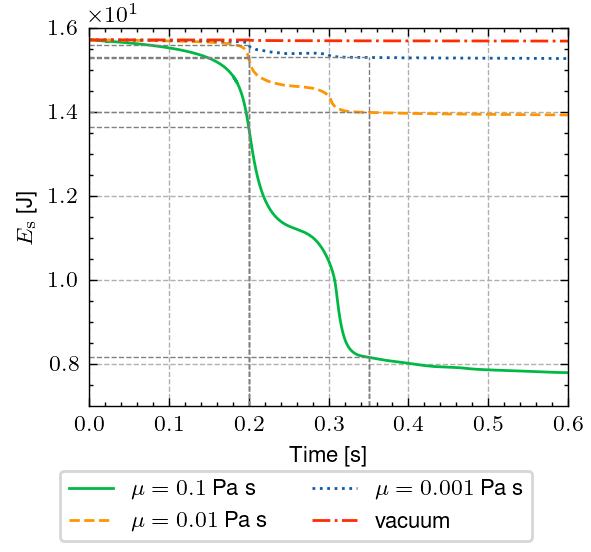}
  \includegraphics[valign=t]{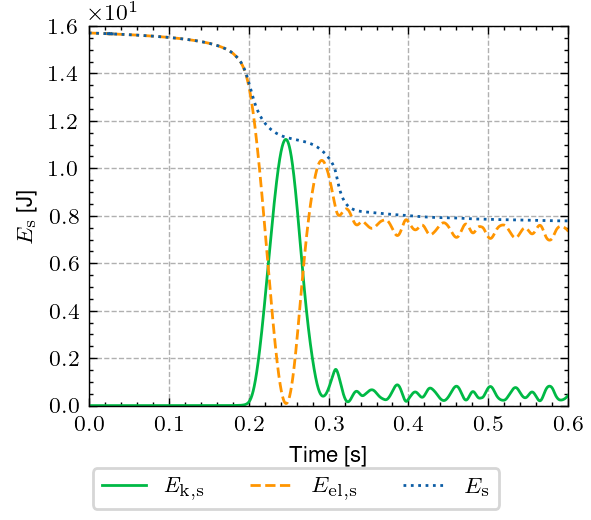}
  \caption{Comparison of the elastic energy $E_s$ for different viscosities (left); comparison of the energies $E_{\rm k, s}$,  $E_{\rm e, s}$ and $E_{\rm s}$ using the viscosity  $\mu=0.1$~Pa\,s (right). All results are obtained using the elastic modulus $G = 50$~kPa, timestep $\Delta t = 10^{-4}$~s, $\meshlev{3}{600}$ and $\eicref$ strategy.}
  \label{fig:energy1}
\end{figure}

Finally, we compare how the ball dissipates the total energy $E_{\rm s}$ with increasing viscosity; see Figure~\ref{fig:energy1}. 
First, during the rebound phase, the energy is pumped from the kinetic energy $E_{\rm k, s}$ to the elastic energy $E_{\rm el, s}$, and then it is transferred back to the kinetic energy.

By comparing the kinetic energy $E_{\rm k,s}$ before the `impact' at $t_{\rm touch}=0.2$\,s, and after the `impact' at $t_{\rm after}=0.35$\,s we can obtain the coefficient of restitution $e$ (ratio between the final and initial velocities, introduced by Newton \cite{Weir2008}) by
\begin{equation}
  e=\sqrt{\frac{E_{\rm k,s}(t_{\rm after})}{E_{\rm k,s}(t_{\rm touch})}}.
\end{equation}

As expected, the coefficient of restitution decreases with increasing viscosity, as we can see in Table~\ref{tab:coefrestitution}. The values correspond to those experimentally observed; see \cite{Goldsmith1960}.
\begin{table}[H]
  \centering
  \begin{tabular}{c|c|c|c}
    $\mu$ [Pa s]   & $E_{\rm k, s}(0.2)$ [J]   & $E_{\rm k, s}(0.35)$ [J]   &  $e$ \\
    \hline
    $0.1  $ &   $13.499$  & $7.543 $ &  $0.748$ \\
    $0.01 $ &   $15.247$  & $13.149$ &  $0.929$ \\
    $0.001$ &   $15.593$  & $14.438$ &  $0.962$ \\
  \end{tabular}
  \caption{Comparison of coefficient of restitution for different viscosities using the elastic modulus $G = 50$~kPa, timestep $\Delta t = 10^{-4}$~s, $\meshlev{3}{600}$ and $\eicref$ strategy.}
  \label{tab:coefrestitution}
\end{table}

\section{Level-set Method}
For the numerical implementation of the fully Eulerian approach of FSI described by \eqref{eulereq1}--\eqref{eulereq2}, we employ the conservative level-set method with reinitialization, see \cite{olsson2005} for details. Similarly, as in \cite{liu2001, gravina2022}, the sharp interface between the fluid and the solid is replaced by the diffuse one, and the method is used to track the interface. In particular, the global unknowns are enhanced by adding the scalar function level-set $l$ distinguishing between the fluid and the solid, i.e.
\begin{equation}
  l(x,t)=\Big\{
  \begin{matrix}
    1,&\textrm{if}\ x\in\Omega_{\rm s},\\
    0,&\textrm{if}\ x\in\Omega_{\rm f}.
  \end{matrix}
\end{equation}
Next, the scalar $l$ is smoothly blurred across the interface with a characteristic thickness $\varepsilon$ and denoted by $l_\varepsilon$. As the solid moves through the fluid, the regularized level-set $l_\varepsilon$ is advected by the velocity, i.e.,
\begin{equation}\label{advection_levelset}
  \pder{l_\varepsilon}{t}+\vecko\cdot\nabla l_\varepsilon = 0.
\end{equation}

\begin{figure}
  \begin{center}
    \includegraphics{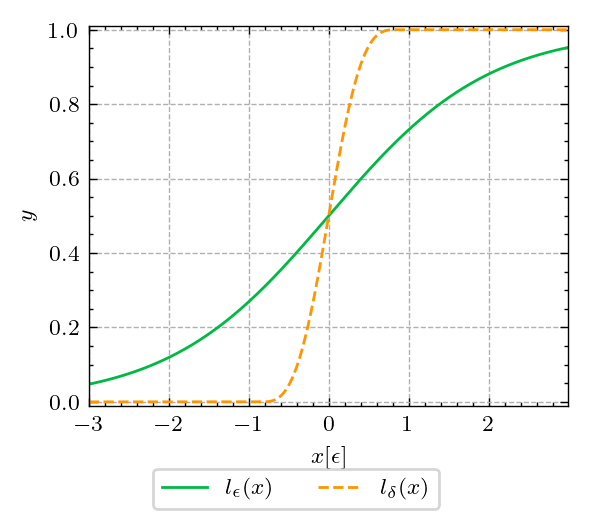}
    \caption{One-dimensional comparison of the original level-set $l_\varepsilon=\left(1+\tanh\left(\frac{x}{2\varepsilon}\right)\right)/2$ and the focused level-set $l_\delta$ given by \eqref{focused-levelset} using $\delta = 0.2$.}
    \label{fig_diffuse_interface}
  \end{center}
\end{figure}

To improve the stability of the method and achieve a good resolution of the interfacial zone, we reinitialize the level-set function during the simulations by solving the following equation
\begin{equation}
  \div\left(\bar{l}_\varepsilon(1-\bar{l}_\varepsilon)\frac{\nabla l_\varepsilon}{|\nabla l_\varepsilon|}\right)-\varepsilon\Delta l_\varepsilon=0,
\end{equation}
where $\bar{l}_\varepsilon$ is the solution of \eqref{advection_levelset}. Next, the solution $l_\varepsilon$ is assigned to $\bar{l}_\varepsilon$, further evolving according to \eqref{advection_levelset}. The thickness of the diffuse interface of such reinitialized level-set corresponds to the parameter $\varepsilon$, see Figure~\ref{fig_diffuse_interface}. 
Originally, the material parameters $\rho,\mu,G$ change sharply across the interface \eqref{differing_material_parameters}, now they change smoothly across the interface, i.e.
\begin{align}
  \rho(l_\varepsilon):=l_\delta\rho_{\rm s}+(1-l_\delta)\rho_{\rm f},\
  \mu(l_\varepsilon):=l_\delta\mu_{\rm s}+(1-l_\delta)\mu_{\rm f},\
  G(l_\varepsilon):=l_\delta G_{\rm s}+(1-l_\delta)G_{\rm f}.
\end{align}
To make the change sharper, we employ the focused level-set $l_\delta$ given by
\begin{linenomath}
  \begin{equation}
    \label{focused-levelset}
    l_{\delta}(l_{\varepsilon}) := \left\{
      \begin{matrix}
        0& \hspace*{13mm}l_\varepsilon < 0.5-\delta,\\
        \frac12 + \frac{l_\varepsilon-0.5}{2\delta} + \frac{1}{2\pi}\sin\left(\frac{\pi(l_\varepsilon-0.5)}{\delta}\right)&0.5-\delta\leq l_\varepsilon \leq 0.5+\delta,\\
        1&\hspace*{13mm}l_\varepsilon > 0.5+\delta.
      \end{matrix}
    \right.
  \end{equation}
\end{linenomath}
In the calculations, the parameter $\delta$ is set to 0.2, for which a very sharp change of materials is obtained, but the numerical scheme is still stable. One-dimensional comparison of the original level-set $l_\varepsilon$ and focused level-set $l_\delta$ is shown in Figure~\ref{fig_diffuse_interface}.

The problem \eqref{eulereq1}--\eqref{eulereq2} is very complex, mainly due to the solution of the evolution equation for the tensor $\mathbb{B}$. To reduce the size of the problem, we apply the following simplifications. Since the flow is quite slow, we omit the convective terms both in the balance of linear momentum \eqref{euler_balance} and the evolution equation for $\mathbb{B}$. The second is potentially problematic and does not allow for reliable stress resolution in the ball during the ``advecting'' phases (before impact and after impact). However, during the short rebound phase, which is part of the process we are primarily interested in, the closest-to-contact part of the ball is static. At this phase, the error caused by omitting the convective time derivative in the evolution of $\mathbb{B}$ in the solid is negligible. 

Next, since the ball deformation is very small, $\mathbb{B}_{\rm s}\sim\II$ and thus $(\nabla\vecko)\BB_{\rm s}\sim\nabla\vecko$ in the solid; and in the fluid, where no elastic stress is present, we prescribe $\mathbb{B}_{\rm f}=\II$. Finally, the global left Cauchy-Green $\mathbb{B}=l_\varepsilon\mathbb{B}_{\rm s}+(1-l_\varepsilon)\mathbb{B}_{\rm f}$ is a mixture of the left Cauchy-Green in the solid $\mathbb{B}_{\rm s}$ and in the fluid $\mathbb{B}_{\rm f}$. In this way, we can explicitly express $\BB$ as a function of the velocity $\vecko$ and level-set $l_\varepsilon$ and insert it into the balance of linear momentum. Thus, we solve the following set of governing equations
\begin{align}
  \label{fsi-mass}
  \div\ \vecko &=0,\\
  \label{fsi-levelset}
  \pder{{l_{\varepsilon}}}{t} + \vecko \cdot \nabla {l_{\varepsilon}} &=0, \qquad \div\left(\bar{l}_{\varepsilon}(1-\bar{l}_{\varepsilon})\frac{\nabla{l_{\varepsilon}}}{|\nabla{l_{\varepsilon}}|}\right) - \varepsilon \Delta \bar{l}_{\varepsilon} = 0,
\end{align}
\begin{align}
  \label{fsi-momentum}
  \rho(l_{\varepsilon}) \pder{\vecko}{t} &= \div\ \mathbb{T},\qquad \mathbb{T} = -p\II + 2 \mu(l_{\varepsilon}) \mathbb{D} + G(l_{\varepsilon}) \mathbb{B}^d,\\
  \label{fsi-B}
  \pder{\mathbb{B}_{\rm s}}{t} & = 2 \mathbb{D},\qquad \mathbb{B}_{\rm f}=\II,\qquad \mathbb{B} = l_{\varepsilon}\mathbb{B}_{\rm s} + (1-l_{\varepsilon})\mathbb{B}_{\rm f}.
\end{align}

The problem is implemented using the FEniCS finite element code \cite{fenics}. The following two subsections describe the space and time discretizations in detail. The final problem \eqref{fsi-mass}--\eqref{fsi-B} deals only with three global unknowns: velocity $\vecko$, pressure $p$, and level-set $l_\varepsilon$. The velocity-pressure pair is approximated by inf-sup stable P2/P1 Taylor-Hood elements, and P2 elements approximate the level-set. The whole problem is solved in a fully monolithic manner. The non-linearities are treated with the Newton scheme that takes advantage of the automatic differentiation available in FEniCS. The resulting set of linear equations is solved with the MUMPS direct solver \cite{mumps}.

\subsection{Space discretization}
The whole domain $\Omega$ is discretized with regular triangles. To capture a sharp interface using a diffuse interface, we must use as fine elements as possible. Earlier, in \cite{gravina2022}, we discretized the domain with an almost uniform mesh. This approach enabled us to describe the rebound of the elastic ball in the Newtonian fluid of the dynamic viscosity 0.2 Pa.s. 
Here, we first mesh the domain uniformly with $125\times 125$ crossed elements. Next, the mesh is refined adaptively close to the interface to be able to compute a lower dynamic viscosity value $\mu$. This is also done by employing the distance function $d(x)$, which is at time $t=0$ initialized and further transported by
\begin{equation}
  \pder{d}{t}+\vecko\cdot\nabla d=0.
\end{equation}
It turns out that for the implementation, it is more efficient and stable to solve the transport in the following way. First, we solve the problem for $\vecko, p, l_{\varepsilon}$ on a given mesh. Next, the velocity $\vecko$ is used to compute the change of (the Lagrangian) displacement $\Delta{\bf u}=\Delta t \ \vecko$, the mesh is moved by this displacement using the FEniCS function {\tt ALE.move()} together with the distance function $d$.

\begin{figure}[!h]
  \centering
  \includegraphics[width=4.95cm]{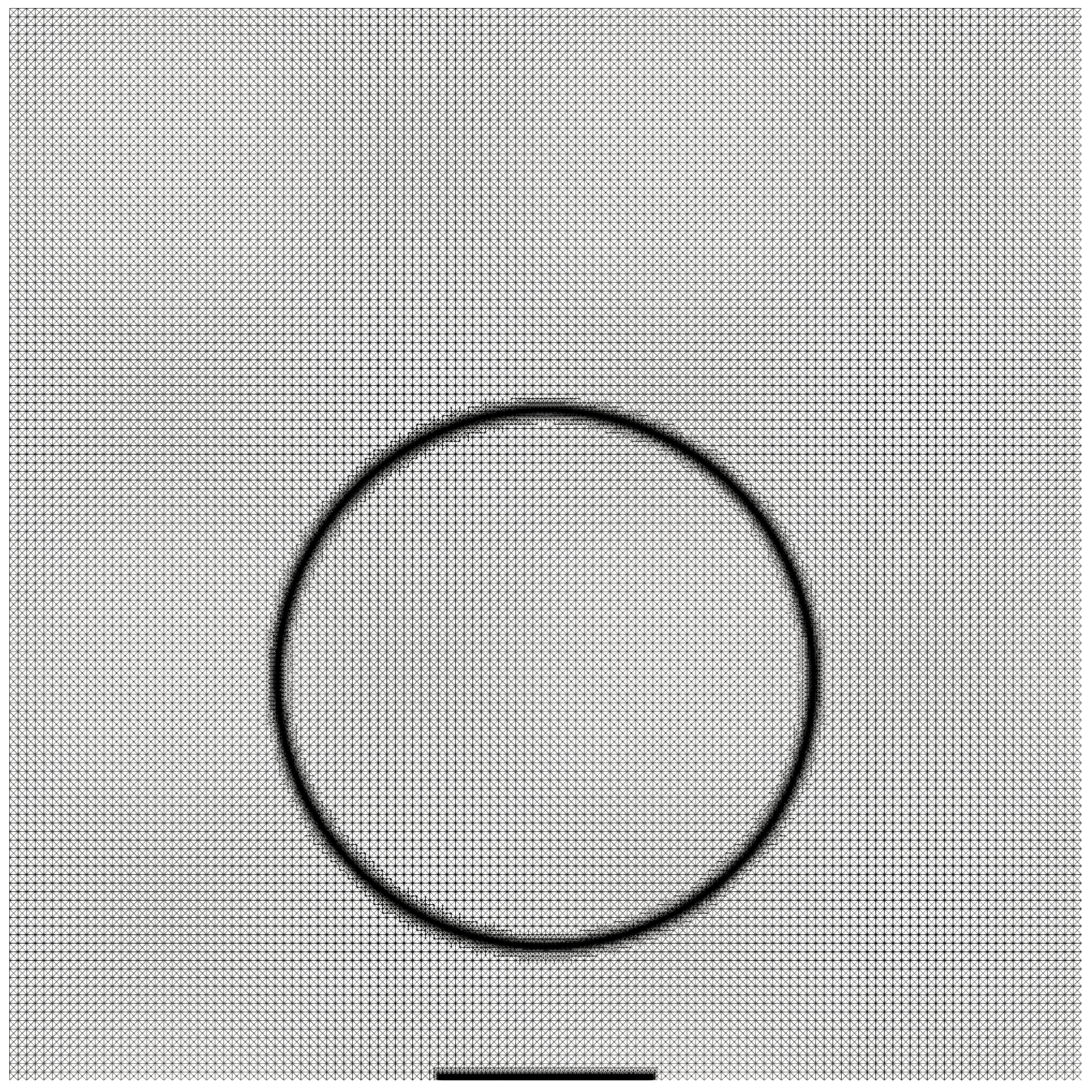}
  \includegraphics[width=4.95cm]{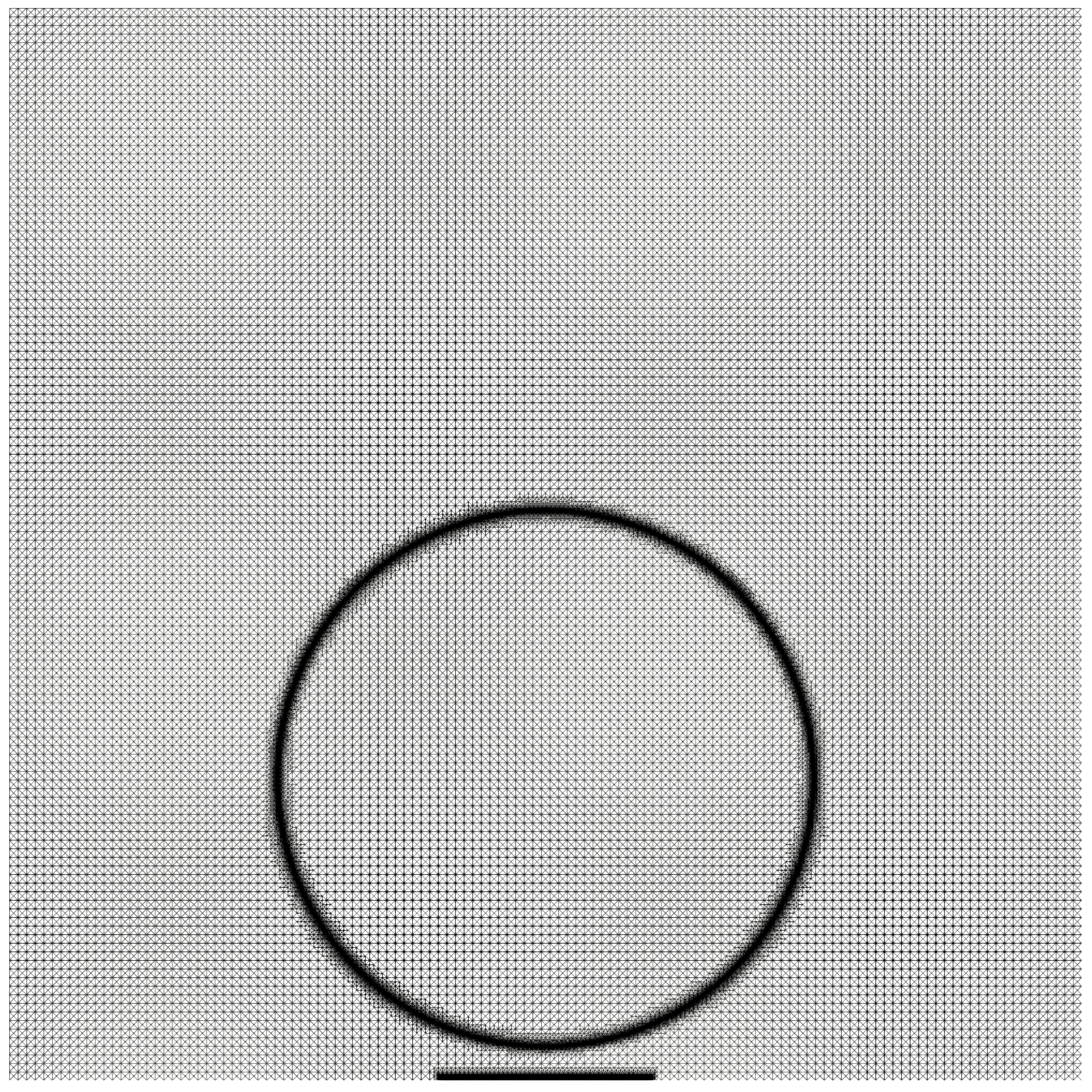}
  \includegraphics[width=4.95cm]{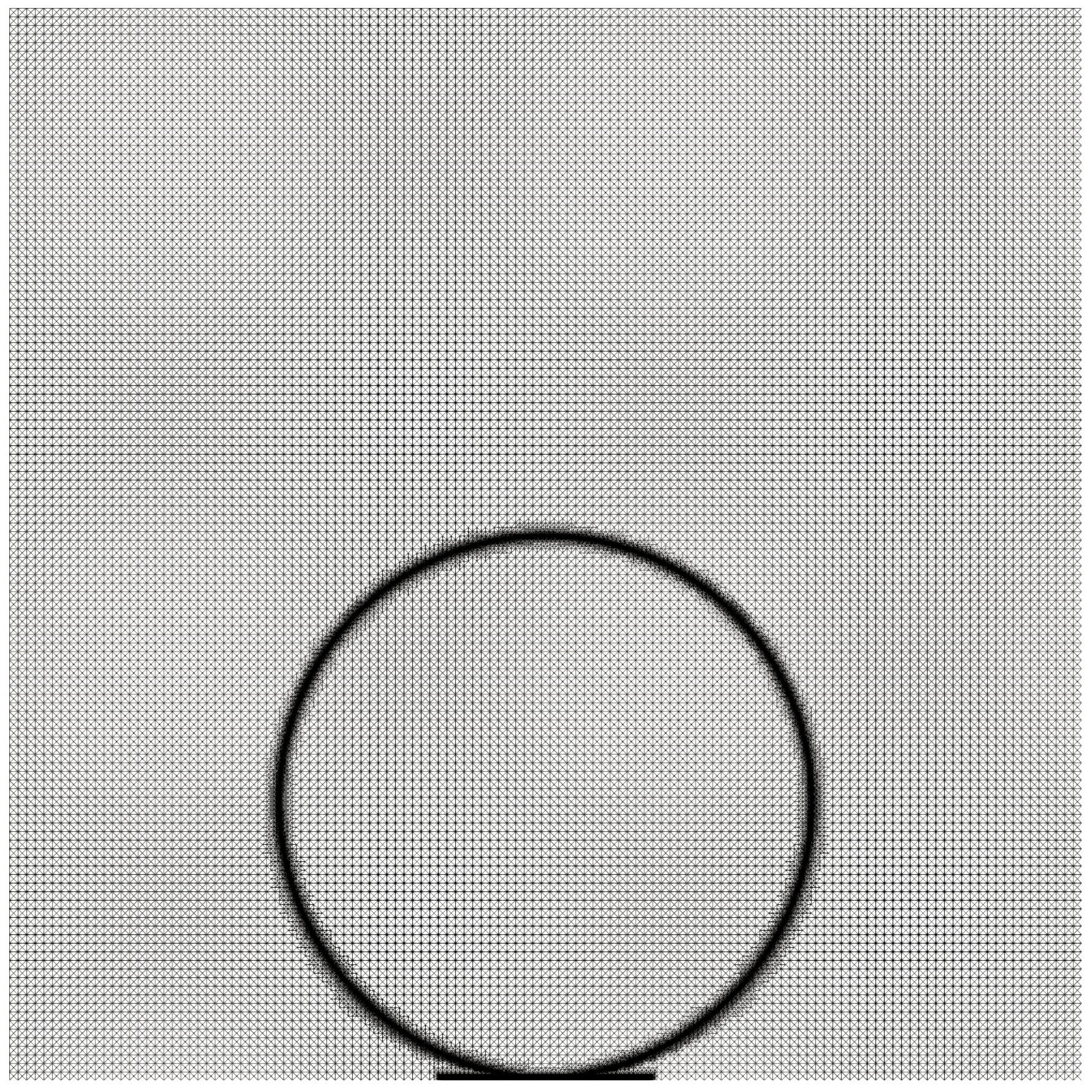}\\
  (a)\hspace*{4.6cm}(b)\hspace*{4.6cm}(c)
  \caption{Remeshing strategy for the fully Eulerian approach: (a) initial mesh, (b) mesh at $t=0.15$\,s, (c) mesh at $t=0.304$~s.}
  \label{fig:mesh_levelset}
\end{figure}

Next, the mesh is five times refined in the vicinity of the interface (thus $h_{\rm min}=10^{-4}$~m) by combining two conditions: 1) the distance function is closer than $1.5\,\varepsilon$, 2) the level-set is close to the value 0.5. Finally, all unknowns $\vecko, p, l_{\varepsilon}$ and $d$ are projected to the new mesh. This re-meshing strategy works very well; see Figure~\ref{fig:mesh_levelset} for the evolution of the mesh.

\medskip

\newpage
\subsection{Time discretization}
The total time interval $[0,T]$ is divided into $N$ subintervals $[t_0,t_1],\dots,[t_{N-1},t_N]$, where $0=t_0, t_1, t_2, \dots, t_N=T$, where $t_n-t_{n-1}=\Delta t=2\times 10^{-4}$~s. 
All time derivatives are approximated with the implicit Euler time scheme (see \eqref{backward_Euler}), and the global left Cauchy-Green tensor at the $n-$th time step is explicitly obtained locally by computing
\begin{equation}
  \mathbb{B}_{\rm s}^{n} = \mathbb{B}_{\rm s}^{n-1} + 2 \mathbb{D}^{n} \Delta t, \quad \mathbb{B}_{\rm f}^{n} = \II, \quad \mathbb{B}^n=l_{\varepsilon}\mathbb{B}_{\rm s}^n + (1-l_{\varepsilon})\mathbb{B}_{\rm f}^n
\end{equation}
and then inserted into the balance of linear momentum (note that $\mathbb{B}^n$ depends on $\DD^n$ at the current time step).

\subsection{Results}
\begin{figure}[!htbp]
  \begin{center}
    \includegraphics[width=7.0cm]{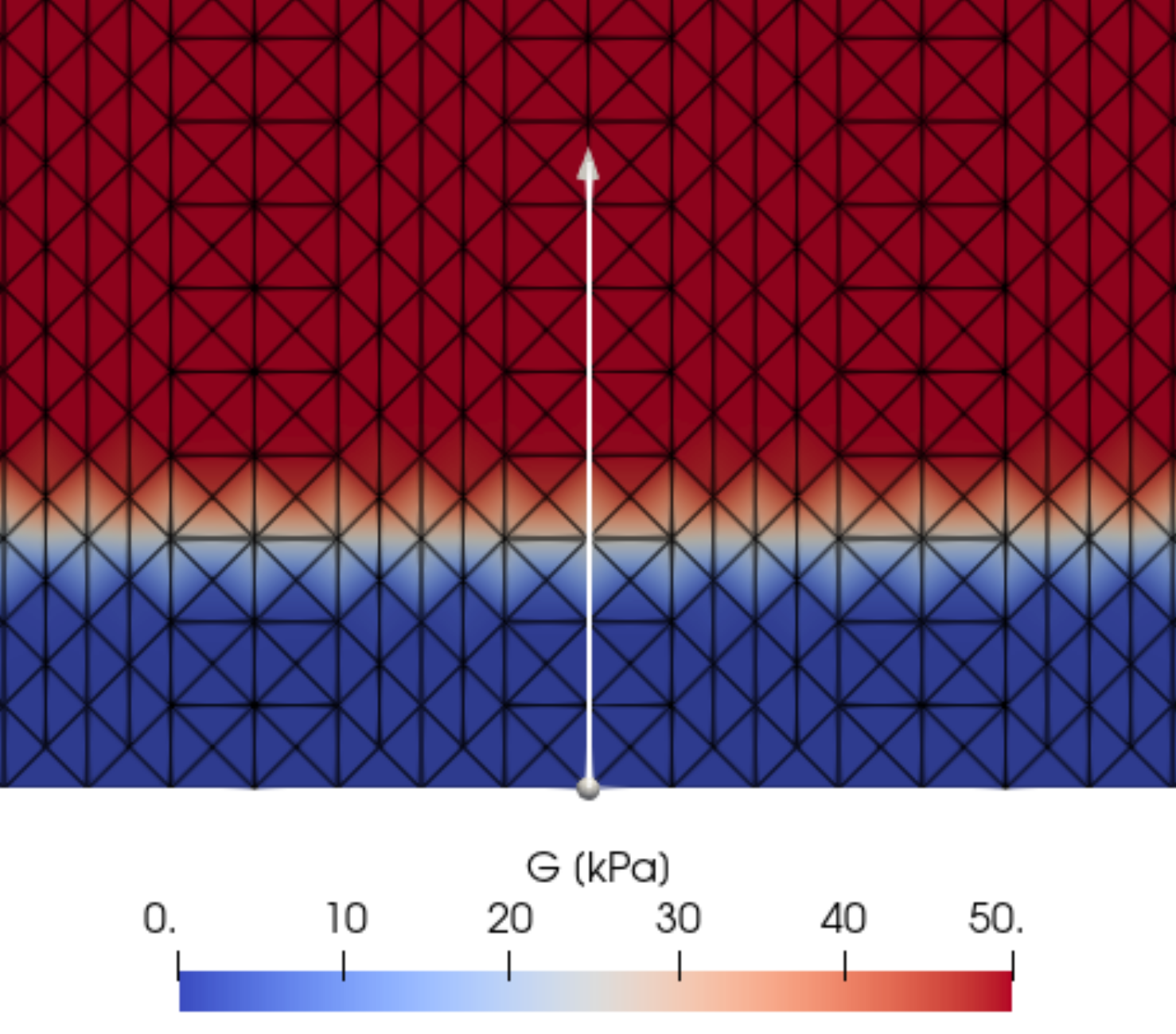}\hspace*{0.5cm}\raisebox{0.0cm}{\includegraphics{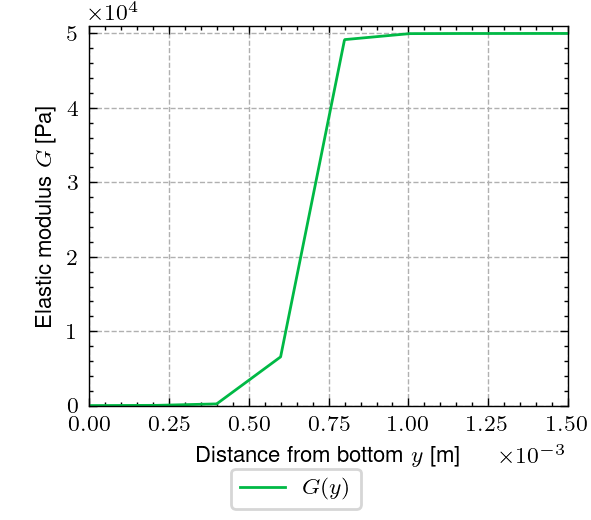}}\\
    (a)\hspace*{6.4cm}(b)
    \caption{(a) Screenshot of the elastic modulus $G$ at the time of rebound $t=0.304$~s. (b) Dependence of the elastic modulus $G$ along the white arrow shown in the screenshot (a).}\label{fig:GalongLine}
  \end{center}
\end{figure}

The results obtained with the level-set method demonstrate a rebound without contact. Indeed, at the time of rebound ($t=0.304$~s), there is a thin layer of pure fluid, see Figure~\ref{fig:GalongLine} that shows the dependence of the elastic modulus $G$ along the white arrow shown in Figure~\ref{fig:GalongLine}a). Clearly, the layer of pure fluid under the ball at the time of rebound is more than three elements thick.

\subsection{Comparison of ALE method and level-set}
\begin{figure}[H]
  \centering
  \includegraphics{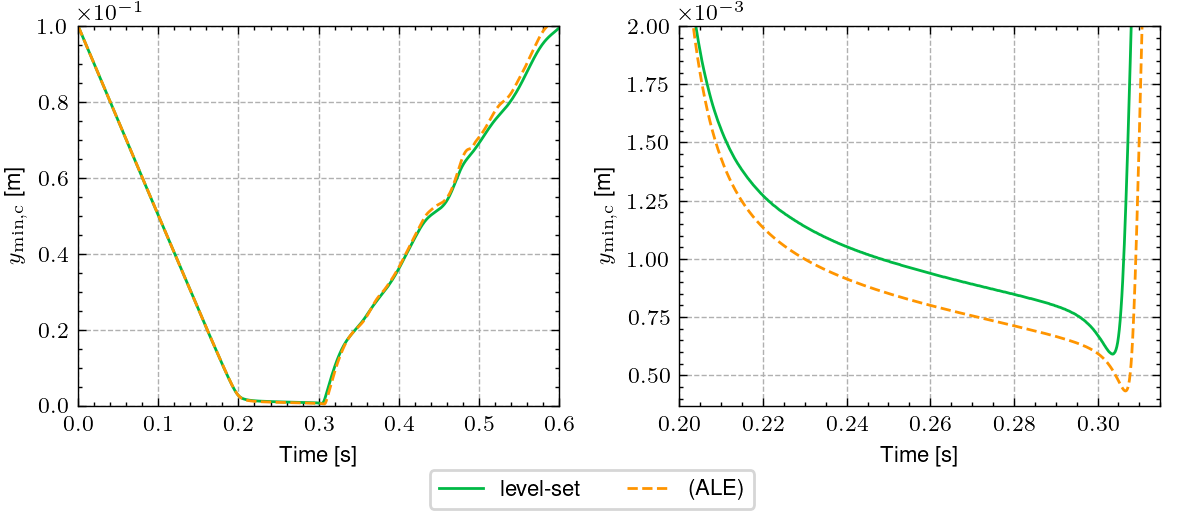}
  \caption{Comparison of time evolution of the $y_{\rm min, c}$ between the level-set and ALE methods. The results are obtained using the viscosity $\mu = 0.1$~Pa\,s the elastic modulus $G = 50$~kPa. The ALE method uses the $\eicref$ strategy on $\meshlev{3}{600}$ with $\Delta t = 10^{-4}$~s.}
  \label{fig:levelset_vs_ale}
\end{figure}

Finally, we compare the results from ALE and the level-set method. We can compare these two approaches in Figure~\ref{fig:levelset_vs_ale}. Due to the non-zero size of the transition between the materials and the slightly different choice of the PDE models, the level-set method does not reach the same distance as the ALE method. Nevertheless, the macroscopic behavior of these two methods is in good agreement. We can see that two methods based on completely different approaches exhibit the rebound with no contact. Thus, we can be confident about the mathematical correctness of the FSI problem. However, the no-contact rebound has not been physically measured yet. It raises the question of whether the no-contact rebound is just a limitation of the continuum mechanics or a real physical phenomenon.

\section{Conclusion}
We have implemented an Updated ALE method with re-meshing strategies using the finite element method. The implementation was successfully tested on the FSI problem of a neo-Hookean solid submerged in an incompressible Newtonian fluid in which an elastic ball is thrown against a rigid wall, and its contactless rebound was studied. The problem was solved for several different fluid viscosities $\mu=0.1, 0.01, 0.001$ Pa\,s.

First, we studied the convergence in space and time for the largest viscosity $\mu=0.1$ Pa\,s. In space, we checked how the solution converges for different (i) approximations of the ball, (ii) global mesh refinement, and (iii) adaptive refinement strategies (geometrical, Eikonal equation). The $\qualityref$ refinement approach ensured that the elements below the ball are regular; however, there is typically only one layer of the elements. On the other hand, the $\eicref$ approach ensured that there was a layer of at least four elements below the ball. We observed that the solution did not change with a finer mesh and used a better $\eicref$ approach. Thus, the solution converged in space. We also observed the convergence in time for the finest mesh and the best (i.e., $\eicref$) refinement strategy. 

Next, we studied how the solution depends on the material parameters such as fluid viscosity $\mu$ and elastic modulus $G$. Due to the Updated ALE method using the $\eicref$ approach, it was possible to decrease the viscosity to $\mu=0.001$ Pa\,s.
As the viscosity decreases, the ball bounces off at a lower distance from the wall (for $\mu=0.001$ Pa\,s the distance is very small, $y_{\rm min}\doteq 41 \mu$m), and we checked that the solution converges to the solution with the ball in a vacuum, for which we implemented the problem using the augmented Lagrangian method.
As the ball is stiffer (elastic modulus $G$ increases), the minimum distance slightly decreases, and the pressure singularity, which is responsible for the contactless rebound, increases.

Finally, we employed another prominent method used to simulate fluid-structure interactions -- the purely Eulerian approach based on the level-set method, which showed very good agreement with the adaptive ALE scheme for $\mu=0.1$ Pa\,s. Smaller viscosities were not achieved with this method.

We conclude that a contactless rebound is certainly an artifact of the PDE. Nevertheless, as the rebound is far from being completely understood, it remains an open question of whether the here considered PDEs suffice to reproduce a real-world situation, even on the macroscopic level. Other aspects as roughness of the materials, could influence the dynamics drastically. However, what our experiments show is that with decreasing viscosity the rebound is clearly converging towards the bounce in vacuum. Consequently the local mechanics that can be read of the numerical experiment is likely to be of qualitative relevance, while its quantitative relevance is still to be investigated.

Since our implementation provided consistent results, the computed problem could potentially be a good setting for a FSI benchmark. We plan to contact other FSI groups, agree on the setting of a benchmark, and compare different numerical approaches. 

%%%
\paragraph{Acknowledgement} J.F. S.S. and K.T.\ have been supported by the ERC-CZ grant LL2105
of the Ministry of Education, Youth, and Sport of the Czech Republic and by Charles University
Research program No. UNCE/SCI/023. J.F. has been supported by Specific University Research project SVV-2023-260711. S.S and K.T. are members of the Necas Center for Mathematical Modeling.
This work has been supported by the Ministry of Education, Youth and Sports of the Czech Republic through the e-INFRA CZ (ID:90140).

\appendix
\section{Rebound of elastic ball in vacuum}\label{Appendix_vacuum}
In this abstract, we provide the details of implementation for the rebound of an elastic ball occupying the domain $\Omega_{\rm s}$ in the reference configuration. As in the FSI problem, the ball is thrown against the wall but is not surrounded by an incompressible fluid. The ball is assumed to be incompressible neo-Hookean that is described by the potential energy functional (cf. \eqref{firstPiola1})
\begin{equation}
  E_{\rm s}[{\bf u}, P]=\int_{\Omega_{\rm s}} \psi_{\rm s}\,{\rm d}X,\quad \psi_{\rm s}=\frac{G}{2}\left(|\mathbb{F}|^2-2\right) + P(\det{\mathbb{F}}-1),
\end{equation}
where $\psi_{\rm s}$ is the potential energy density, $G$ is the elastic modulus, $\mathbb{F}=\mathbb{I}+\nabla{\bf u}$ is the deformation gradient and pressure $P$ is the Lagrange multiplier responsible for the incompressibility.

The evolutionary problem complemented with the initial condition (zero displacement and prescribed initial velocity) can be formulated in the variational form
\begin{equation}
  \int_{\Omega_{\rm s}}\rho_R\displaystyle\pder{^2{\bf u}}{t^2}\cdot {\bf \varphi}_{\bf u}\,{\rm d}X +
  \int_{\Omega_{\rm s}}\displaystyle\pder{\psi_{\rm s}}{\nabla{\bf u}}\cdot \nabla{\bf \varphi}_{\bf u}\,{\rm d}X + \int_{\Omega_{\rm s}}(\det\mathbb{F}-1)\,{\varphi}_{P}\,{\rm d}X=0
\end{equation}
that holds for all suitable test functions ${\bf \varphi}_{\bf u}$ and $\varphi_P$ and equivalently, it can be written using the Gateaux derivative of the functional $E_{\rm s}$, i.e.
\begin{equation}
  \int_{\Omega_{\rm s}}\rho_R\displaystyle\pder{^2{\bf u}}{t^2}\cdot {\bf \varphi}_{\bf u}\,{\rm d}X +
  \delta_{{\bf u}}E_{\rm s}+\delta_{P}E_{\rm s}=0.
\end{equation}
The two Gateaux derivatives represent the necessary conditions for the saddle point problem
\begin{equation}\label{saddlepoint1}
  \min_{\bf u} \max_{P} E_{\rm s}.
\end{equation}
So far, we have not considered the wall. This can be easily formulated by saying that no point of the ball can occur under the wall. Since the wall is horizontal (satisfying $y=0$), we require that the gap $g$ (in the current configuration) between the ball and the wall is non-negative, i.e.
\begin{equation}\label{constraint}
g = Y + u_y \geq 0,
\end{equation}
where $Y$ is the position in the reference configuration and $u_y$ is the $y$ component of the displacement. To deal with this inequality constraint, we employ the augmented Lagrangian method; see \cite{Bertsekas2014} for more details on this method and \cite{AlartCurnier,StupkiewiczLengiewicz} for applying this method on contact. 

The saddle point problem \eqref{saddlepoint1} with inequality constraint \eqref{constraint} is transformed into the following smooth and unconstrained double saddle-point problem
\begin{equation}
  \min_{\bf u} \max_{P} \max_{\lambda} \mathcal{L}[{\bf u}, P, \lambda],\quad \mathcal{L}[{\bf u}, P, \lambda]=E_{\rm s}[{\bf u}, P] + \int_{\Omega_{\rm s}} l({\bf u}, \lambda)\,{\rm d}X,
\end{equation}
where $l({\bf u}, \lambda)$ denotes a continuously differentiable function responsible for the exact fulfillment of the inequality constraint \eqref{constraint} and $\lambda$ is the corresponding Lagrange multiplier. Function $l({\bf u}, \lambda)$ is then defined in the following way
\begin{equation}
  l({\bf u}, \lambda)=\Bigg\{
  \begin{matrix}
    \left(\lambda+\dfrac{1}{2}\rho g\right)g,&\lambda+\rho g\leq 0,\\[7pt]
    -\dfrac{\lambda^2}{2\rho},&\lambda+\rho g>0,
  \end{matrix}
\end{equation}
where $\rho>0$ is a regularization parameter. Finally, the problem considered is formulated in the variational form
\begin{equation}
  \int_{\Omega_{\rm s}}\rho_R\displaystyle\pder{^2{\bf u}}{t^2}\cdot {\bf \varphi}_{\bf u}\,{\rm d}X +
  \delta_{{\bf u}}\mathcal{L}+\delta_{P}\mathcal{L}+\delta_{\lambda}\mathcal{L}=0
\end{equation}
that is equivalent to
\begin{equation}
  \int_{\Omega_{\rm s}}\rho_R\displaystyle\pder{^2{\bf u}}{t^2}\cdot {\bf \varphi}_{\bf u}\,{\rm d}X +
  \int_{\Omega_{\rm s}}\displaystyle\pder{\psi_{\rm s}}{\nabla{\bf u}}\cdot \nabla{\bf \varphi}_{\bf u}\,{\rm d}X + \int_{\Omega_{\rm s}}(\det\mathbb{F}-1)\,{\varphi}_{P}\,{\rm d}X
  + \int_{\Omega_{\rm s}}\pder{l}{\bf u}\cdot{\varphi}_{\bf u}\,{\rm d}X
  + \int_{\Omega_{\rm s}}\pder{l}{\lambda}\,{\varphi}_{\lambda}\,{\rm d}X=0
\end{equation}
for all suitable test functions $\varphi_{\bf u}$, $\varphi_P$, $\varphi_{\lambda}$.

The variational form is implemented in FEniCS finite element code using regular triangles. The second time derivative of the displacement ${\bf u}$ is approximated by using an additional unknown velocity 
\begin{equation}
  {\bf v}=\displaystyle\pder{\bf u}{t}.
\end{equation}
The obtained set of first-order equations in time is approximated using the Glowinski time scheme; see Subsection \ref{Sec:Glowinski}.

The displacement ${\bf u}$, velocity ${\bf v}$, and Lagrange multiplier $\lambda$ are approximated with piecewise quadratic elements, and pressure $P$ is approximated with piecewise linear elements. FEniCS provides the automatic differentiation that is used to compute the derivatives. The piecewise function $l$ is implemented using the FEniCS function {\tt conditional}. The Newton method treats the non-linearities, and the consequent set of linear equations is solved with the direct solver MUMPS.

\bibliographystyle{elsarticle-harv}
\bibliography{bibliography}

\end{document}